\documentstyle[12pt,amstex,amssymb]{amsart}

\newtheorem{thm}[subsection]{Theorem}
\newtheorem{prop}[subsection]{Proposition}
\newtheorem{lemma}[subsection]{Lemma}
\newtheorem{cor}[subsection]{Corollary}

\def\L{{\cal L}}
\def\N{{\cal N}}
\def\C{{\frak C}^{\scriptstyle{\text{nil}}}}
\def\g{{\frak g}}

\def\t{{\frak t}}
\def\N{{\cal N}}

\def\z{{\frak z}}

\def\Lie{\mathop{\fam0 Lie}}

\def\Mat{\mathop{\fam0 Mat}\nolimits}
\def\sl{\mathop{\frak sl}\nolimits}

\def\SO{\mathop{\text{SO}}\nolimits}
\def\SL{\mathop{\text{SL}}\nolimits}

\def\GL{\mathop{\text{GL}}\nolimits}
\def\subtitle#1. {{\medskip\bf#1\par\nobreak\smallskip}}
\def\proclaim#1. {\medbreak\bgroup\noindent\bf#1. \it}
\def\endproclaim{\egroup
\ifdim\lastskip<\medskipamount\removelastskip\medskip\fi}
\newcount
=0
\def\citedef#1 {\advance\citation by1
 \expandafter\edef\csname#1\endcsname{{\the\citation}}
 \checkendcitedef}
\def\checkendcitedef#1{\ifx#1\endcitedef\else\citedef#1\fi}
\def\cite#1{\csname#1\endcsname}
\citedef Bas Bar BR Bo Bou Br Ca  F1 F2 Ge Gra Gr HSp I Jac Jan J Kr Le
LS Mc Na Po1 Po2 P86 P89 P1
P2 P02 Rich1 Rich2  Sh Spa Spr SpSt St SFB
\endcitedef
\newtoks\nextauth
\newif\iffirstauth
\def\checkendauth#1{\ifx\endauth#1
        \iffirstauth\the\nextauth
        \else{} and \the\nextauth\fi,
    \else\iffirstauth\the\nextauth\firstauthfalse
        \else, \the\nextauth\fi
        \expandafter\auth\expandafter#1\fi}
\def\auth#1 #2 {\nextauth={#1 #2}\checkendauth}
\newif\ifinbook
\newif\ifbookref
\def\nextref#1 {\bookreffalse\inbookfalse
    \bibitem[\cite{#1}]{}
    \firstauthtrue
    \ignorespaces}
\def\paper#1{{\it#1,}}
\def\In#1{\inbooktrue In #1,}
\def\book#1{\bookreftrue{\it#1,}}
\def\journal#1{#1\ifinbook,\fi}
\def\bookseries#1{#1,}
\def\Vol#1{\ifbookref Vol. #1,\else\ifinbook Vol. #1,\else{\bf#1}\fi\fi
    \space\ignorespaces}

\def\publisher#1{#1,}
\def\Year#1{\ifbookref #1.\else\ifinbook #1,\else(#1)\fi\fi
    \space\ignorespaces}
\def\Pages#1{\ifinbook pp. #1.\else #1.\fi}

\begin{document}

\title{Nilpotent commuting varieties of reductive\\ Lie algebras}

\vspace{0.7cm}
\author{Alexander Premet}
\thanks{\nonumber{\it Mathematics Subject Classification} (1985 {\it revision}).
Primary 20G05}
\address{Department of Mathematics, University of Manchester, Oxford Road,
M13 9PL, UK} \email{sashap@@ma.man.ac.uk}
\begin{abstract}
\noindent Let $G$ be a connected reductive algebraic group over an
algebraically closed field $k$ of characteristic $p\ge 0$, and
$\g=\text{Lie}\,G$. In positive characteristic, suppose in
addition that $p$ is good for $G$ and the derived subgroup of $G$
is simply connected. Let $\N=\N(\g)$ denote the nilpotent variety
of $\g$, and
$\C(\g)\,:=\,\{(x,y)\in\N\times\N\,\vert\,\,[x,y]=0\},$ the
nilpotent commuting variety of $\g$. Our main goal in this paper
is to show that the  variety $\C(\g)$ is equidimensional. In
characteristic $0$, this confirms a conjecture of Vladimir
Baranovsky; see [\cite{Bar}]. When applied to $ \GL(n),$ our
result in conjunction with an observation in [\cite{Bar}] shows
that the punctual (local) Hilbert scheme ${\cal
H}_{n}\subset\text{Hilb}^{n} ({\mathbb P}^2)$ is irreducible over
any algebraically closed field.
\end{abstract}
\maketitle

\section{\bf Introduction}
Let $k$ be an algebraically closed field of characteristic $p\ge
0$. The purpose of this note is to confirm Baranovsky's conjecture
[\cite{Bar}, p.~4] which states that all irreducible components of
the nilpotent commuting variety $\C(\g)$ of a complex semisimple
Lie algebra $\g$  have the same dimension, equal to $\dim\,\g$,
and are parametrised by the distinguished nilpotent orbits in
$\g$. The main result in [\cite{Bar}] confirmed the conjecture for
$\g\,=\,\sl(n)$ (an earlier proof of the irreducibility of
$\C(\sl(n))$ in [\cite{Gr}] was incomplete). Notably, it is
observed in [\cite{Bar}] that in any characteristic the number of
irreducible components of $\C(\sl(n))$ equals the number of
irreducible components of the punctual Hilbert scheme ${\cal
H}_n\subset\text{Hilb}^n ({\mathbb P}^2)$; see [\cite{Bar},
Remark~1].

The punctual
Hilbert scheme ${\cal H}_n$ parametrises the ideals of colength $n$
in the ring of formal power series
$k[[x,y]]$. In characteristic $0$, the scheme ${\cal H}_n$
is known to be irreducible
for more than 25 years thanks to the work of Brian\c{c}on [\cite{Br}].
This was extended to the case where $p>n$ by
Iarrobino [\cite{I}];
see also [\cite{Gra}]. These results enabled Baranovsky to deduce that
$\C(\sl(n))$ is irreducible for $p=0$ and $p>n$. It should be mentioned here that
very recently a more direct proof of the irreducibility
of $\C(\sl(n))$ was found by Basili in [\cite{Bas}]. It works for
$p\ge n/2$ and $p=0$ implying the irredicibility of
${\cal H}_n$ for $p$ in that range.

In this note we give a direct proof of
Baranovsky's conjecture entirely in the framework of Lie Theory.
In view of
[\cite{Bar}, Remark~1], this will enable us to conclude that
${\cal H}_n$ is irreducible over any algebraically closed
field. For $p<n/2$, this appears to be a new result in Algebraic Geometry.
Our approach
also provides a much shorter and more elementary proof of the
irreducibility of ${\cal H}_n$ over $\mathbb C$.
In principle, it can be used for
investigating the connected
components of ${\cal H}_n$ over other locally compact fields; see
[\cite{I}] and [\cite{Bar}, Remark~2].

Let $G$ be a connected reductive algebraic group over $k$ and
$\g=\text{Lie}\,G$. We assume throughout the paper that the
derived subgroup $(G,G)$ is simply connected and $p$ is good  for
the root system of $G$. When $p>0$, the Lie algebra $\g$ carries a
natural $(\text{Ad}\,G)$-equivariant $[p]$-mapping $x\mapsto
x^{[p]}$. In this case, an element $e\in\g$ is called nilpotent if
$e^{[p]^N}=0$ for $N$ large enough. When $p=0$, we say that $e$ is
nilpotent if $e\in[\g,\g]$ and the endomorphism $\text{ad}\,e$ is
nilpotent. The variety of all nilpotent elements in $\g$ is
denoted by $\N$. It is well-known that $\N$ is an irreducible
Zariski closed subset in $\g$ of dimension $\dim\,G-\text{rk}\,G$
and $G$ acts on $\N$ with finitely many orbits; see [\cite{Rich1}].
Moreover, the Bala--Carter theory holds in good
characteristic and the $(\text{Ad}\,G)$-orbits in $\N$ are
described in the same way as over $\mathbb C$. In other words, any
nilpotent element in $\g$ is $(\text{Ad}\,G)$-conjugate to a
Richardson element in a distinguished parabolic subalgebra of
$\text{Lie}\,L$ where $L$ is a Levi subgroup of $G$; see
[\cite{Ca}, \cite{Po1}, \cite{Po2}, \cite{P02}]. Let
$$\C(\g)\,:=\,\{ (x,y)\in\N\times\N\,|\,\,
[x,y]=0\}\subset\g\times\g,$$ the nilpotent commuting variety of
$\g$. Obviously, the Zariski closed set $\C(\g)$ is preserved by
the diagonal action of $G$ on $\g\times\g$.

Given a closed subgroup $H$ in $G$ we denote by $Z_{H}(e)$ the
centraliser of $e$ in $H$. As usual we denote by $R_u(H)$ the unipotent
radical of $H$. An element $e\in\N$ is called
{\it distinguished} if the connected component of $Z_{(G,G)}(e)$ is
unipotent, that is if $\Lie\,Z_{(G,G)}(e)\subset\N$. According to
the main result of the Bala--Carter theory, any distinguished
nilpotent element in $\g$ is Richardson in a distinguished
parabolic subalgebra of $\g$. Note that for any distinguished
$e\in\N$ we have $(e,\,\Lie\, Z_{(G,G)}(e))\subset\C(\g)$. We denote
by ${\frak C}(e)$ the Zariski closure of $G\cdot (e,\,\Lie\, Z_{(G,G)}(e))$
in $\C(\g)$. Our main result in this paper is the following:

\smallskip

\noindent {\bf Theorem.} {\it Let} $e_1,\ldots,e_r$ {\it be
representatives of the distinguished nilpotent orbits in} $\g$.
{\it The closed sets}\, ${\frak C}(e_1),\ldots,{\frak C}(e_r)$
{\it are pairwise distinct and all have the same dimension equal to}
$\dim\,(G,G)$. {\it Moreover,}
$\C(\g)\,=\,{\frak C}(e_1)\cup\ldots\cup {\frak C}(e_r).$

\medskip

A well-known result of Richardson [\cite{Rich2}] asserts that for
$p=0$ the  whole commuting variety ${\frak C}(\g)\,:=\,\{(x,y)\in
\g\times \g\,|\,\,[x,y]=0\}$ of $\g$ coincides with the Zariski
closure of $G\cdot (\Lie\,T \times\Lie\,T)$ in $\g\times\g$ where
$T$ is a maximal torus in $G$. As a consequence, this variety  is
always irreducible. The case where $p$ is good was recently settled  by
P. Levy in his PhD thesis; see also [\cite{Le}]. An important
long-standing conjecture asserts that the variety ${\frak C}(\g)$
is always normal and the ideal defining ${\frak C}(\g)$ is given
by the ``obvious'' quadratic equations. It is not hard to see that
$\C(\g)$ coincides with the special fibre of the quotient map
${\frak C}(\g)\,\rightarrow\, {\frak C}(\g)/\!\!/G$. For
$p=0$, a Chevalley Restriction Theorem holds for ${\frak C}(\g)$
[\cite{J}]; it says that the affine variety ${\frak C}(\g)/\!\!/G$
is isomorphic to $(\text{Lie}\,T \times\text{Lie}\,T)/\!\!/W$ where the
action of the Weyl group $W=N_{G}(T)/T$ on $\Lie \,T\times
\Lie\,T$ is induced by the diagonal action of $G$ on
$\g\times\g$. It would be useful to have an analogue of this in
positive characteristic.
\smallskip

To emphasize the elementary nature of our proof we first consider
the case where $p=0$ (this is done in Section~2). In the modular
case our argument goes along the same lines but is technically
much more involved (see Section~3). This is mainly due to
inseparability problems and a rather unusual behaviour of the
orbit map:  in small characteristic there exist nontrivial
commuting ${\frak sl}_2$-triples $(e_1,h_1, f_1)$ and
$(e_2,h_2,f_2)$ in $\g$ such that $e_1$ and $e_1+e_2$ are
conjugate under the adjoint action of $G$. The latter is, in our
opinion, the main reason why the irredicibility of ${\cal H}_n$ is
harder to establish for $p$ small. To tackle this problem we  go
case-by-case and look closely at the centralisers of nilpotent
elements. For exceptional types, we have to use the computations
in [\cite{P1}].

As a motivation for further study, we
mention that the nilpotent commuting variety and its higher analogues
play a very important r{\^o}le in the cohomology
theory of the Frobenius kernels of $G$. It is proved in [\cite{SFB}] that
$\C(\g)$ is homeomorphic to the spectrum of maximal
ideals of the Yoneda algebra $\bigoplus_{i\ge 0}\,H^{2i}(G_2,k)$  of
the second Frobenius kernel $G_2$ of $G$
provided that $p$ is sufficiently large.  The variety  $\C(\g)$ is also important
for the study of support varieties of modules over
reduced enveloping algebras of $\g$; see [\cite{P2}].

In Section~4, we prove that the punctual Hilbert scheme ${\cal H}_n$
is irreducible over any algebraically closed field; see Corollary~4.1.
We also show, in (4.2),  that the nilpotent variety of any finite dimensional
restricted Lie algebra is equidimensinal. This result is then used in (4.3) to
estimate the dimension of $\C(\g)$ in the case where $p$ is a bad prime for
$G$.

\smallskip

\noindent {\bf Acknowledgement.} I would like to thank
V.~Ginzburg, J.C.~Jantzen, D.~Panyushev, D.~Timashev and
{\`E}.B.~Vinberg for their interest, comments and informations. I
am also grateful to P.~ Levy from whom I learned that $\GL(2)$
acts on the commuting variety of $\g$.

\section{\bf The classical case}

\noindent {\bf 2.1.} Unless otherwise indicated we assume in this
section that $p=0$. In this case, $\g\,=\,[\g,\g]\oplus\z$ where
$\z$ is the Lie algebra of the connected centre of $G$.  Thus no
generality will be lost by assuming that $G$ is semisimple. Then
$G$ and $\g$ are both defined over $\mathbb Q$ and hence so is
$\C(\g)$. Therefore, all irredicible components of $\C(\g)$ are
defined over the field of algebraic numbers. Thus it can be
assumed in what follows that $k=\mathbb C$.

Let $e$ be a nilpotent element in $\g$.
Let $\z(e)$ denote the centraliser of $e$ in $\g$.
By the
Jacobson--Morozov theorem, $e$ can be embedded into an $\sl_2$-triple $(e,h,f)$ in
$\g$. By the $\sl_2$-theory, all eigenvalues of the semisimple
endomorphism $\text{ad}\,h$ are integers. For $i\in\mathbb Z$ we let
$\g(i;h)\,=\,\{x\in\g\,|\,\,[h,x]=ix\}$. Then $\g\,=\,\bigoplus_{i\in\mathbb Z}\,\g(i;h)$
and $[\g(i;h),\g(j;h)]\subseteq\g(i+j;h)$ for all $i,j\in\mathbb Z$.

Since all nonzero elements in $\z(e;i)\,:=\,\g(i;h)\cap\z(e)$ are
heighest weight vectors for ${\frak s}\,:=\,{\mathbb C}e\oplus
{\mathbb C}h\oplus {\mathbb C}f$ we have the inclusion
$\z(e)\subset \bigoplus_{i\ge 0}\,\g(i;h)$ (again, by the
$\sl_2$-theory).  The Lie algebra $\z(e;0)$ is nothing but the
centraliser of $\frak s$ in $\g$ hence reductive. Therefore, the
ideal $\bigoplus_{i>0}\,\z(e;i)$ of $\z(e)$ coincides with the
nilradical of $\z(e)$. As a consequence, the variety $$\N\cap
\z(e)\,=\,\N(\z(e;0))+\,\text{nil}\,\z(e) \,\cong\,
\N(\z(e;0))\times\text{nil}\,\z(e)$$ is irreducible. We denote by
${\frak C}(e)$  the Zariski closure of $G\cdot(e,\,\N\cap\z(e))$
in $\g\times \g$.  It is immediate from the definition that
${\frak C}(e)\subseteq \C(\g)$ and ${\frak C}(e)\,=\,{\frak
C}((\text{Ad}\,g)\,e)$ for any $g\in G$. The preceding remark
shows that each ${\frak C}(e)$ is irreducible.

\smallskip

\noindent
{\bf Definition.} We call  a nilpotent element $e\in\g$  {\it almost distinguished} if
the centraliser  of $e$ in $\g$ is a solvable Lie algebra.

\smallskip

Note that $e$ is almost distinguished if and only if the connected
component of the centraliser $Z_G(e)$ is solvable. Any
distinguished nilpotent element in $\g$ is therefore almost
distinguished. The converse, however, holds only for  simple
algebraic groups of small rank. Looking through the tables in
[\cite{Ca}, Chap.~13] one finds out that the exceptional groups
for which the converse also holds are the groups of types
$\mathrm{G}_2$ and $\mathrm{F}_4$. At the other extreme, the Lie
algebra $\sl(n)$ has only one distinguished nilpotent orbit while
there is a bijection between the almost distinguished nilpotent
orbits in $\sl(n)$ and the partitions of $n$ with pairwise
distinct parts.

Since $\z(e;0)$ is reductive, the Lie algebra $\z(e)$ is solvable if and only if
$\z(e;0)$ coincides with the Lie algebra of a maximal torus in $Z_G(e)$.
From this it is
immediate that $e\in\N$ is almost distinguished if and only if
$\N\cap\z(e)\,=\,\bigoplus_{i>0}\,\z(e;i).$
\begin{prop} Each irreducible component of
${\frak C}^{\mathrm{nil}}(\g)$
is of the form ${\frak C}(e)$ for some
almost distinguished $e\in\g$.
\end{prop}
\pf (1) The group $\GL(2)$ acts on $\g\times\g$  via
$$\textstyle{\left({\alpha\atop\gamma}\, {\beta\atop\delta}\right)\cdot (x,y)\,=
\,(\alpha x+\beta y,\gamma x+\delta y)}.$$
As any linear combination of two
commuting elements in $\N$ is again in $\N$, the
variety $\C(\g)$ is $\GL(2)$-invariant. As $\GL(2)$ is a connected group, it fixes each
irreducible component of $\C(\g)$.
In particular, each irreducible component of
$\C(\g)$ is invariant under the involution $\sigma\colon\,(x,y)\mapsto (y,x)$ on $\g\times \g$.

\smallskip

\noindent
(2) Let $e_1,\ldots e_s$ be representatives of the nilpotent orbits in $\g$. Since each
${\frak C}(e_i)$ is $G$-stable we have that
$\C(\g)\,=\,{\frak C}(e_1)\cup\ldots\cup{\frak C}(e_s).$ This implies that each irreducible
component of $\C(\g)$ has the form ${\frak C}(e_i)$ for some $i\le s$. Let $e\in \N$ be
be such that ${\frak C}(e)$ is a component of $\C(\g)$. By part~1, the set ${\frak C}(e)$ is
then $\sigma$-stable.

Let ${\cal O}$ denote the $G$-orbit of $e$. The map
$\pi\colon\,(x,y)\mapsto x$ takes $G\cdot (e,\,\N\cap\z(e))$ onto $\cal O$, hence
${\frak C}(e)\,=\,\overline{G\cdot (e,\,\N\cap\z(e))}$
onto the Zariski closure $\overline{\cal O}\subset \g$. This shows that
$$\N\cap\z(e)=(\pi\circ\sigma)(e,\,\N\cap\z(e))\subset\overline{\cal O}.$$

\noindent
(3) We need to show that the element $e$ is almost distinguished. So
suppose the contrary.
Then $\z(e;0)$ contains a nonzero nilpotent element, say $e_0$.
Note that $e+e_0\in\N\cap\z(e)$ so that $e+e_0\in\overline{\cal O}$, by part~2.
Since $\z(e;0)$ is reductive,
 $e_0$  can be included into an $\sl_2$-triple $(e_0,h_0,f_0)$ in $\z(e;0)$, again by the
Jacobson--Morozov theorem. Let ${\frak s}_0\,=\,{\mathbb
C}e_0\oplus {\mathbb C}h_0\oplus {\mathbb C}f_0$. Since ${\frak
s}_0\subseteq\z(e;0)$, the Lie subalgebras ${\frak s}$ and ${\frak
s}_0$ commute with each other. This enables us to deduce that
$(e+e_0,h+h_0,f+f_0)$ is an $\sl_2$-triple in $\g$. Applying the
automorphisms $\exp\,(\lambda\,\text{ad}\,(h+h_0))\circ
\exp\,(\lambda^{-1}\,\text{ad}\,h)$ from the adjoint group of $\g$
to $e+e_0$ we deduce that $e+{\mathbb C}^*e_0\subset{\cal
O}(e+e_0)$, the $G$-orbit of $e+e_0$. Then $e\in\overline{{\cal
O}(e+e_0)}$. As a result, $\overline{{\cal
O}(e+e_0)}\,=\,\overline{\cal O}$ showing that $e$ and $e+e_0$ are
$G$-conjugate. But then the semisimple elements $h$ and $h+h_0$
are $G$-conjugate too; see [\cite{Ca}, Prop.~5.6.4] for example.

\smallskip

\noindent
(4) Since $\g(0;h)$ is reductive we have
that $\g(0;h)\,=\,[\g(0;h),\g(0;h)]\oplus\z'$ where $\z'$ is the centre of  $\g(0;h)$.
The subspaces
 $\z'$ and $[\g(0;h),\g(0;h)]$ are orthogonal to each other with respect to the Killing form
$\kappa$ of $\g$. Since ${\frak s}_0\subseteq [\g(0;h),\g(0;h)]$ and $h\in\z'$ it must be that
$\kappa(h,h_0)=0$. On the other hand, $\kappa(h,h)>0$ and $\kappa(h_0,h_0)>0$
as all eigenvalues of $\text{ad}\,h$ and $\text{ad}\,h_0$ are in $\mathbb Z$.
But then $$\kappa(h+h_0,h+h_0)=\kappa(h,h)+\kappa(h_0,h_0)>\kappa(h,h).$$
Since $h$ and $h+h_0$ are $G$-conjugate (see part~3) this is impossible.
By contradiction, the proposition follows.
\qed

\smallskip

\noindent {\bf 2.2.} In a sense, our quest requires that we attach
to an arbitrary $\sl_2$-triple $(e,h,f)$ in $\g$ a nice
$\sl_2$-triple $(\tilde{e},\tilde{h},\tilde{f})$ with $\tilde{e}$
distinguished. This will be achieved with the help of the
Bala--Carter theory; see [\cite{Ca}, Chap.~5] and [\cite{P02}]. In
this subsection, we work with an arbitrary reductive group $G$
over $k$ assuming only that $p=\text{char}\,k$ is good for $G$.

Fix a maximal torus $T$ in $G$ and let $\Phi$ denote the root
system of $G$ relative to $T$. Let
$\Pi\,=\,\{\alpha_1,\ldots,\alpha_l\}$ be a basis of simple roots
in $\Phi$, $\Phi^+$ be the positive system in $\Phi$ associated
with $\Pi$, and $X_{*}(T)$ be the lattice of cocharacters of $T$
(the latter contains all coroots $\alpha^\vee$). For a subset
$I\subseteq \{1,\ldots,l\}$ let $L_I$ (resp., $P_I$) stand for the
standard Levi (respectively, parabolic) subgroup of $G$
corresponding to $I$. Let $\Phi_I$ be the root system of $L_I$
relative to $T$. The set $\Pi_I\,:=\,\{\alpha_i\,|\, i\in I\}$ is
then a basis of simple roots in $\Phi_I$.

Given two subsets $I\supseteq J$ in $\{1,\ldots,l\}$ we denote by $P_{I,J}$ the
standard parabolic subgroup of $L_I$ associated with $J$. Let ${\frak l}_I\,=\,
\Lie\,L_I$, ${\frak p}_I\,=\,\Lie\,P_I$,
${\frak p}_{I,J}\,=\,\Lie\,P_{I,J}$, and ${\frak u}_{I,J}\,=\,\Lie\,R_u(P_{I,J})$.
Note that $P_{I,J}\,=\,L_{J}\cdot R_u(P_{I,J})$ is a Levi decomposition in $P_{I,J}$.
According to [\cite{Ca}, Prop.~5.8.2],
for each  pair $(I,J)$ with $\{1,\ldots,l\}\supseteq I\supseteq J$ the inequality
\begin{eqnarray}
\dim\,(L_I,L_I)\cap L_J\,\ge\,\dim\,R_u(P_{I,J})/\big(R_u(P_{I,J}),R_u(P_{I,J})\big)
\end{eqnarray}
holds. A subgroup $P_{I,J}$ is said to be {\it distinguished} in $L_I$
if $$\dim\,(L_I,L_I)\cap
L_J\,=\,\dim\,R_u(P_{I,J})/\big(R_u(P_{I,J}),R_u(P_{I,J})\big).$$
Let ${\cal P}(\Pi)$ be the set of all pairs $(I,J)$
with$\{1,\ldots,l\}\supseteq I \supseteq J$ such that $P_{I,J}$ is
distinguished in $L_I$. For $(I,J)\in{\cal P}(\Pi)$ we denote by
${\cal O}(I,J)$ the nilpotent orbit in $\g$ containing a
Richardson element in ${\frak u}_{I,J}$. The main result of the
Bala-Carter theory (as extended to the present case in
[\cite{Po1}, \cite{Po2}, \cite{P02}]) asserts that
$\N\,=\,\bigcup_{(I,J)\in{\cal P}(\Pi)}{\cal O}(I,J)$. Thus we may
assume in what follows that $e$ is a Richardson element in ${\frak
u}_{I,J}$.

Given $\mu\in X_*(T)$ and  a $\mu(k^*)$-invariant subalgebra
$\frak m$ in $\g$ we denote by ${\frak m}(i;\mu)$ the subspace of
all $x\in\frak m$ such that $(\text{Ad}\,\mu(t))\,x=t^{i}x$ for
all $t\in k^*$ (here $i\in\mathbb Z$). As explained in
[\cite{P02}, Sect.~2] (for example) there exists a cocharacter
$\lambda_{I,J}\in \sum_{i\in I}{\mathbb Z}\alpha_i^\vee$ such that
$$\alpha_i(\lambda_{I,J}(t))\,\,=\,\,\left\{\begin{array}{ll}
1&\mbox{ if $\,i\in J$},\\ t^2&\mbox{ if $\,i\in I\setminus J$}
\end{array}
\right.$$ for all $t\in k^*$.
When $I=\{1,\ldots,l\}$ (that is when $P_J$ is distinguished in $G$)
we denote the cocharacter $\lambda_{I,J}$ by $\lambda_J$.
Since the orbit $(\text{Ad}\,P_{I,J})\,e$ meets
${\frak l}_{I}(2;\lambda_{I,J})$, by [\cite{Ca}, Prop.~5.8.5], we may (and will) assume
that $e\in {\frak l}_{I}(2;\lambda_{I,J})$.

Our next result is implicit in [\cite{LS}]. We give a direct proof
for the benefit of the reader.
\begin{lemma} [cf.
$\mathrm{[\cite{LS}, Prop.~1.12]}$]
For each $(I,J)\in{\cal P}(\Pi)$ the  parabolic subgroup
$P_{J}$ is distinguished in $G$.
\end{lemma}
\pf
Following [\cite{Ca}, (5.8)] define a function $\eta_J\colon\,\Phi\rightarrow
2\mathbb Z$ by $$\eta_J(\alpha_i)\,\,=\,\,\left\{\begin{array}{ll}
0&\mbox{ if $\,i\in J$},\\ 2&\mbox{ if $\,i\in\{1,\ldots,l\} \setminus J$}
\end{array}
\right.$$ and extending to arbitrary roots by linearity.
Applying [\cite{Ca}, Prop.~5.8.1] to $P_J\subset G$ and $P_{I,J}\subset L_I$
and using the fact that $(I,J)\in{\cal P}(\Pi)$
we deduce that
\begin{eqnarray*}
\dim\,R_u(P_{J})/\big((R_u(P_{J}),R_u(P_{J})\big)&=&
\text{Card}\,\{\alpha\in\Phi\,|\,\eta_{J}(\alpha)=2\}\\
&\ge &|\Pi\setminus \Pi_I|+\text{Card}\,\{\alpha\in\Phi_I\,|\,\eta_{J}(\alpha)=2\}\\
&=&l-|I|+\dim\,R_u(P_{I,J})/\big((R_u(P_{I,J}),R_u(P_{I,J})\big)\\
&=&l-|I|+\dim\,(L_I,L_I)\cap L_J.
\end{eqnarray*}
On the other hand, $\dim\,(G,G)\cap L_J=l+|\Phi_I|=
l-|I|+\dim\,(L_I,L_I)\cap L_J.$ Hence
$$\dim\,R_u(P_{J})/\big((R_u(P_{J}),R_u(P_{J})\big)\ge \dim\,(G,G)\cap L_J.$$
Applying (1) with $I=\{1,\ldots,l\}$ yields that $P_J$ is distinguished in $G$
as desired.
\qed

\smallskip

\noindent
{\bf 2.3.} From now  we assume in this section that $k=\mathbb C$ and $\g$ is
semisimple.
Let $(I,J)\in{\cal P}(\Pi)$ and let $e\in{\frak l}_{I}(2;\lambda_{I,J})$ be
such that the orbit
$(\text{Ad}\,P_{I,J})\,e$ is open in ${\frak u}_{I,J}$.
It is immediate from [\cite{Ca}, Cor.~5.2.4] that the map
$$\text{ad}\,e\colon\,{\frak l}_{I}(-2;\lambda_{I,J})\,
\longrightarrow\,[{\frak l}_{I},{\frak l}_{I}](0;\lambda_{I,J})$$
is bijective. This implies that there are $h\in \Lie\,\lambda_{I,J}({\mathbb C}^*)$
and $f\in{\frak l}_{I}(-2;\lambda_{I,J})$ such that $(e,h,f)$ is an $\sl_2$-triple in $\g$.
Moreover, $\alpha_i(h)=0$ for $i\in J$ and $\alpha_i(h)=2$ for $i\in I\setminus J$.

By Lemma~2.2, $P_J$ is a distinguished parabolic subgroup of $G$.
Then $\lambda_J\in X_*(T)$ is well-defined and ${\frak
p}_J\,=\,\bigoplus_{i\ge 0}\,\g(2i;\lambda_J)$. By [\cite{Ca},
Prop.~5.8.5], the subspace $\g(2;\lambda_J)$ contains a Richardson
element of ${\frak p}_J$. We pick such an element and call it
$\tilde{e}$. There is $\tilde{h}\in\Lie\,\lambda_J({\mathbb C}^*)$
such that $\alpha_i(\tilde{h})=0$ for $i\in J$ and
$\alpha_i(\tilde{h})=2$ for $i\in\{1,\ldots,l\}\setminus J$. Since
the map
$\text{ad}\,\tilde{e}\colon\,\g(-2;\tilde{h})\,\longrightarrow\,
\g(0;\tilde{h})$ is bijective,  by [\cite{Ca}, Cor.~5.2.4], there
is $\tilde{f}\in\g(-2;\tilde{h})$ such that
$(\tilde{e},\tilde{h},\tilde{f})$ is an $\sl_2$-triple in $\g$.

By construction,
$\g(i;\lambda_J)\,=\,\g(i;\tilde{h})$ and $\g(i;\tilde{h})\cap
{\frak l}_{I}\,=\,{\frak l}_{I}(i;h)$ for all $i\in\mathbb Z$. From this it follows that
\begin{eqnarray} [h,\tilde{h}]=0,\quad\,\,
[\tilde{h},e]=2e,\quad\,\
\g(0;\tilde{h})\,=\,\g(0;\lambda_J)\,=\,{\frak l}_J.
\end{eqnarray}
So the centraliser $\z(e)$ is $(\text{ad}\,\tilde{h})$-invariant.
From (2.1) we know that  $\text{nil}\,\z(e)$ coincides with
$\bigoplus_{i>0}\,\g(i;h)\cap\z(e)$, hence is also preserved by
$\text{ad}\,\tilde{h}$.
\begin{lemma}
The endomorphism ${\mathrm ad}\,\tilde{h}$ acts invertibly on
the nilradical ${\mathrm nil}\,\z(e)$.
\end{lemma}
\pf
Indeed, it follows from (2) that \begin{eqnarray*}
\text{nil}\,\z(e)\cap\text{Ker}\,\text{ad}\,\tilde{h}&=&
\text{nil}\,\z(e)\cap \g(0,\tilde{h}) \,=\,
\text{nil}\,\z(e)\cap {\frak l}_J\\&=&\text{nil}\,\z(e)\cap{\frak l}_I(0;h)
\subseteq \text{nil}\,\big({\frak l}_I(0;h)\cap\z(e)\big). \end{eqnarray*}
Since the Lie algebra ${\frak l}_I(0;h)\cap\z(e)$ is reductive, by our discussion
in (2.1), the RHS must be zero. The result follows.
\qed

\smallskip

\noindent
{\bf Remark.} Using the terminology of the theory of induced orbits
[\cite{LS}] one can say that $\tilde{e}$ is a nice correspondent of $e$ in
the $G$-orbit
$\text{Ind}_{L_I,\,P_I}({\cal O}_I(e))$ where ${\cal O}_I(e)=
(\text{Ad}\,L_I)\,e$.

\smallskip

\noindent
{\bf 2.4.} Let $V$ be a finite dimensional vector space over $\mathbb C$.
Given an ideal $I$ of the symmetric algebra $S(V^*)$
we denote by $\text{gr}\,I$ the homogeneous ideal of $S(V^*)$ with the property
that $g\in\text{gr}\,I\cap S^r(V^*)$ if and only if there is $\tilde{g}\in I$ such that
$\tilde{g}-g\in\bigoplus_{i<r}\,S^i(V^*)$. As usual,
we identify $S(V^*)$ with the algebra of polynomial  functions on $V$.
We denote by ${\cal Z}(\varphi_1,\ldots, \varphi_N)$  the subset in $V$
consisting of all common zeros of $\varphi_1,\ldots, \varphi_N\in S(V^*)$.
Given
a Zariski closed set $X\subseteq V$ we let $I_X$ be the ideal in $S(V^*)$
consisting of all polynomial functions vanishing on $X$ and define
$${\mathbb K}(X)\,:=\,\{v\in V\,|\,f(v)=0\quad\forall\,f\in\text{gr}\,I_X\}.$$
The Zariski closed conical set ${\mathbb K}(X)\subseteq V$ is known as the
{\it associated cone} to $X$. The proof of our next result will rely on
a few elementary properties of the correspondence $X\,\longmapsto\,{\mathbb K}(X)$.
Our reference here is [\cite{Kr}, Chap.~1, (4.2)].

Given an $\sl_2$-triple $(e,h,f)$ in $\g$ we denote by ${\cal
S}(h,e)$ the Zariski closure of the orbit $G\cdot(h,e)$ in
$\g\times\g$. Identify $\g\times\g$ with $\g\oplus\g$. The next
result will play a crucial  r{\^o}le in our proof of Baranovsky's
conjecture.
\begin{prop}
Let $(e,h,f)$ be an $\sl_2$-triple in $\g$ with  $e$  distinguished.
Then there exist pairwise non-conjugate
distinguished nilpotent elements
$e_1,\ldots,e_k\in\g$ such that
${\mathbb K}\,\big({\cal S}(h,e)\big)\,=\,{\frak C}(e_1)\cup\ldots\cup {\frak C}(e_k).$
\end{prop}
\pf (1) As $[h,e]=2e$, any pair $(u,v)\in G\cdot(h,e)$ has the property that
$[u,v]=2v$. Therefore, ${\cal S}(h,e)\subset\{(x,y)\in\g\times\g\,|\,[x,y]-2y=0\}.$
In view of [\cite{Kr}, Chap.~2, Sect.~4.2, Theorem~2(b)] this entails that
\begin{eqnarray}
{\mathbb K}\big({\cal S}(h,e)\big)\subset{\mathbb K}\big(
\{(x,y)\in\g\times\g\,|\,\,[x,y]-2y=0\}\big)\,\subseteq \,{\frak C}(\g).
\end{eqnarray}
Let $x_1,\ldots,x_n$ and $y_1,\ldots,y_n$ be coordinate functions
on $\g\times\{0\}$ and $\{0\}\times\g$, respectively. Let
$f_1,\ldots,f_l$ be free homogeneous generators of the invariant
algebra $S(\g^*)^G$. It is well-known that $\N\,=\,{\cal
Z}(f_1,\ldots,f_l)$. As ${\cal S}(h,e)\subset
(\text{Ad}\,G)\,h\times \N$, all polynomial functions
$f_{i}(x_1,\ldots,x_n)-f_i(h)$ and $f_i(y_1,\ldots,y_n)$ on
$\g\times\g$ vanish on ${\cal S}(h,e)$. This gives ${\mathbb
K}\big({\cal S}(h,e)\big)\subset\N\times\N$. Combining this with
(3) one obtains \begin{eqnarray}{\mathbb K}\big({\cal
S}(h,e)\big)\,\subseteq\, \C(\g).
\end{eqnarray}

\smallskip

\noindent (2) By the definition of ${\cal S}(h,e)$, the orbit
morphism $$\psi\colon\,G\,\longrightarrow\,{\cal S}(h,e),\quad\ \
\, \psi(g)\,= \, \big((\text{Ad}\,g)\,h,(\text{Ad}\,g)\,e\big),$$
is dominant. The fibre $\psi^{-1}(\psi(1))$ is nothing but the
stabiliser of $(h,e)$ in $G$, a closed subgroup of $G$. Its Lie
algebra consists of all $x\in\g$ satisfying $[x,h]= [x,e]=0$,
hence coincides with $\z(e;0)$. From (2.1) we know that  $\z(e;0)$
is reductive. Because $e$ is distinguished, we have
$\z(e)\subset\N$; so it must be that $\z(e;0)=\{0\}$. Then
$\psi^{-1}(\psi(1))$ is finite, implying $\dim\,{\cal
S}(h,e)\,=\,n\,=\,\dim\,\g.$ Since ${\cal S}(h,e)$ is irreducible
we now can apply [\cite{Kr}, Chap.~2, Sect.~4.2, Theorem~2(c)] to
deduce that all irreducible components of ${\mathbb K}\big({\cal
S}(h,e)\big)$ have dimension $n$.

\smallskip

\noindent (3) Let $e_1,\ldots,e_r$ be representatives of the
distinguished nilpotent orbits in $\g$. Let ${\cal O}(e_i)$ denote
the nilpotent orbit in $\g$ containing $e_i$ where $1\le i\le r$.
Since $\pi({\frak C}(e_i))\,=\,\overline{{\cal O}(e_i)}$ for all
$i$, the Zariski closed irreducible sets ${\frak
C}(e_1),\ldots,{\frak C}(e_r)$ are pairwise distinct. By the
definition of ${\frak C}(e_i)$, the morphism $$\xi\colon\,G\times
\z(e_i)\,\longrightarrow\, {\frak C}(e_i),\quad\ \ \,
\xi(g,x)\,=\, \big((\text{Ad}\,g)\,e_i,(\text{Ad}\,g)\,x\big),$$
is dominant. For every $x\in\z(e_i)$ the fibre
$\xi^{-1}(\xi(1,x))$ is nothing but the set of all pairs
$\big(g,(\text{Ad}\,g)^{-1}\,x\big)$ with $g\in Z_{G}(e_i)$. It
follows that $\xi^{-1}(\xi(1))\,\cong\, Z_G(e_i)$ as varieties.
The theorem on the dimension of the fibres of a morphism now gives
$$\dim\,{\frak
C}(e_i)\,=\,\dim\,G+\dim\,\z(e_i)-\dim\,Z_G(e_i)\,=\,n\quad\ \ \
(1 \le i\le r).$$ According to [\cite{Bar}, Theorem~2],
$\dim\C(\g)=n$ and the $n$-dimensional irreducible components of
$\C(\g)$ are parametrised by the distinguished nilpotent conjugacy
classes in $\g$. In conjunction with the above discussion this
yields that ${\frak C}(e_1),\ldots,{\frak C}(e_r)$ are exactly the
$n$-dimensional irreducible components of $\C(\g)$. Combining (4)
with our final remark in part~2 we now conclude, after renumbering
the $e_i$'s if necessary, that ${\mathbb K}\big({\cal
S}(h,e)\big)\,=\,{\frak C}(e_1)\cup\ldots\cup{\frak C}(e_k)$ for
some $k\le r$. \qed

\noindent

\smallskip

\noindent {\bf 2.5.} Now all our tools are in place and we are
ready for the main result of this section.
\begin{thm} Let $\g$ be a complex semisimple Lie algebra of dimension $n$.
All irreducible components of the nilpotent commuting variety
${\frak C}^{\mathrm nil}(\g)$ are $n$-dimensional and the number of components equals
the number of distinguished nilpotent conjugacy classes in $\g$.
More precisely, ${\frak C}^{\mathrm nil}(\g)\,=\,{\frak C}(e_1)\cup\ldots\cup{\frak
C}(e_r)$ where $e_1,\ldots,e_r$ are representatives of the
distinguished nilpotent orbits in $\g$.
\end{thm}
\pf Let $e$ be an almost distinguished nilpotent element in $\g$.
In view of Proposition~2.1 it suffices to show that ${\frak
C}(e)\subseteq {\frak C}(e_i)$ for some $i\le r$. Since ${\frak
C}(e)\,=\,{\frak C}((\text{Ad}\,g)\,e)$ for any $g\in G$ no
generality will be lost by assuming that $e$ satisfies the
conditions of (2.3). We include $e$ into an $\sl_2$-triple
$(e,h,f)$ according to the recipe in (2.3) and then consider the
corresponding $\sl_2$-triple $(\tilde{e},\tilde{h},\tilde{f})$
with $\tilde{e}$ distinguished (see (2.3) for more detail). Thanks
to Proposition~2.4 we are then reduced to show that ${\frak
C}(e)\subseteq {\mathbb K}\big({\cal
S}(\tilde{h},\tilde{e})\big)$.

Let ${\cal S}\,=\,{\cal S}(\tilde{h},\tilde{e})$ and
$L\,=\,Z_{G}(\tilde{h})$. As
$\dim\,\g(0;\tilde{h})\,=\,\dim\,\g(2;\tilde{h})$ and
$\z(\tilde{e};0)\,=\,\{0\}$, the orbit $(\text{Ad}\,L)\,\tilde{e}$
is open in $\g(2;\tilde{h})$. Since $\cal S$ is $L$-stable and
Zariski closed, we then have
$\big(\tilde{h},\,\g(2;\tilde{h})\big)\,\subset\,\cal S.$ Since
$\g(2;\tilde{h})\cap{\frak l}_I\, =\,{\frak l}_I(2;h)$, by (2.3),
we obtain $(\tilde{h},{\mathbb C}e)\subset\cal S$.

By construction, $\tilde{h}\in\Lie\,T\subseteq{\g}(0;h)$, see
(2.3), while from (2.1) we know that
$$\Lie\,R_u(Z_{G}(e))\,=\,\text{nil}\,\z(e)\,=\,\textstyle{\bigoplus}_{i>0}\,\,\z(e;i).$$
As the group $R_u(Z_G(e))$ is unipotent, the orbit $
(\text{Ad}\,R_u(Z_G(e)))\,\tilde{h}$ is Zariski closed in
$\tilde{h}+\text{nil}\,\z(e)$. As $[\text{nil}\,\z(e),\tilde{h}]\,
=\,\text{nil}\,\z(e)$, by Lemma~2.3, it is Zariski open in
$\tilde{h}+\text{nil}\,\z(e)$ too, and hence is the whole of
$\tilde{h}+\text{nil}\,\z(e)$. Applying the operators from
$R_u(Z_G(e))$ to $(\tilde{h},{\mathbb C}e)\subset\tilde S$ we now
derive that $(\tilde{h}+\text{nil}\,\z(e),\,{\mathbb C}
e)\,\subset\,\cal S$. But then $${\mathbb
K}\big((\tilde{h}+\text{nil}\,\z(e),\,{\mathbb C}
e)\big)\,=\,(\text{nil}\,\z(e),\,{\mathbb C}e)\,\subset\, {\mathbb
K}(\cal S),$$ by [\cite{Kr}, Chap.~2, Sect.~4.2, Theorem~2(b)].
Since ${\mathbb K}(\cal S)$ is $G$-stable and $e$ is almost
distinguished in $\g$, our discussion in (2.1) yields
$G\cdot\,(\N\cap\z(e),\,e)\,\subseteq {\mathbb K}(\cal S)$. Hence
${\frak C}(e)\,\subseteq\, {\mathbb K}(\cal S)$
and our proof is complete. \qed

\section{\bf The modular case}

\noindent
{\bf 3.1.} In this section we assume that $p=\text{char}\,k$ is  good for $G$.
Recall that  the derived subgroup $(G,G)$ is simply connected.
So there exist simple, simply connected algebraic $k$-groups
$G_1,\ldots,G_m$ such that $(G,G)\cong G_1\times\cdots\times G_m$.
Let $\g_i=\text{Lie}\,G_i$ where $1\le i\le m$, and $\g'=\text{Lie}\,(G,G)$.
Then $\g'=\g_1\oplus\cdots\oplus\g_m$, a direct sum of restricted Lie algebras.
It is well-known that $\N\subset\g'$; see [\cite{P02}, (2.3)] for example.
From this it is immediate that $$\C(\g)\,=\,\C(\g')\,\cong\,\C(\g_1)\times\cdots\times \C(\g_m).$$
Clearly, a nilpotent element $x=x_1+\cdots+x_m$ with $x_i\in \g_i$ is distinguished in $\g$ if
and only if each $x_i$ is distinguished in $\g_i$. Also,
${\frak C}(x)\,\cong\,{\frak C}_{G_1}(x_1)\times\cdots\times {\frak C}_{G_m}(x_m)$
where ${\frak C}_{G_i}(x_i)$ denotes the Zariski closure of
$G_i\cdot (x_i,\text{Lie}\,Z_{G_i}(x_i))$ in $\g_i\times \g_i$. This observation reduces
computing the irreducible components of $\C(\g)$ to the case where $(G,G)$ is a simple
algebraic group.

Until the end of this section we will thus assume that $G$ is
either $\GL(n)$ or a simple algebraic group of type different from
${\mathrm A}$. Let $X_*(G)$ denote the set of all $1$-parameter
subgroups $\mu\colon k^*\rightarrow G$. Given $\mu\in X_*(G)$ and
$i\in\mathbb Z$ we denote by $\g(i;\mu)$ the subspace of all
$x\in\g$ such that $(\text{Ad}\,\mu(t))\,x=t^i x$ for all $t\in
k^*$. Then $\g\,=\,\bigoplus_{i\in\mathbb Z}\,\g(i;\mu)$, that is
each $\mu\in X_*(G)$ induces  a $\mathbb Z$-grading of the
restricted Lie algebra $\g$. We denote by $Z(\mu)$ the centraliser
of $\mu$ in $G$, a Levi subgroup in $G$.

Let $e$ be a nilpotent element in $\g$. As in Section~2, we let
$\z(e)$ be the centraliser of $e$ in $\g$, ${\frak C}(e)$ be the
Zariski closure of $G\cdot(e,\N\cap\z(e))$, and  ${\cal O}(e)$ be
the $G$-orbit of $e$. Recall from (2.2) that ${\cal
O}(e)\,=\,{\cal O}(I,J)$ for some $(I,J)\in{\cal P}(\Pi)$. More
specifically, there is $g\in G$ such that $(\text{Ad}\, g)\,e$ is
a Richardson element of ${\frak p}_{I,J}$ contained in ${\frak
l}_I(2;\lambda_{I,J})$. Set $\lambda_e:=g\,\lambda_{I,J}\,g^{-1}$,
an element in $X_*(G)$. Since $e\in\g(2;\lambda_e)$, the torus
$\lambda_e(k^*)$ acts on $\z(e)$. For $i\in\mathbb Z$, set
$\z(e;i)\,:=\,\z(e)\cap \g(i;\lambda_e)$. According to
[\cite{Jan}, Sect.~5] and [\cite{P02}, Sect.~2], the group
$Z(\lambda_e)\cap Z_G(e)$ is reductive,
$\text{Lie}\big(Z(\lambda_e)\cap Z_G(e)\big)=\, \z(e;0),$ and
$$\z(e)=\text{Lie}\,Z_G(e)= \,\textstyle{\bigoplus_{i\ge
0}}\,\z(e;i),\qquad\, \text{Lie}\,R_u(Z_G(e))=
\,\textstyle{\bigoplus_{i>0}}\,\z(\g,i).$$ This implies that the
nilpotent variety $\N(\z(e;0))\,=\,\N\cap\z(e;0)$ of $\z(e;0)$ is
irreducible. Hence the varieties $\N\cap\z(e) \,\cong
\,\N(\z(e;0))\times\,{\textstyle{\bigoplus}_{i>0}}\,\z(e;i)$ and
${\frak C}(e)$ are irreducible, too.

\begin{lemma} Let $e\in\N$ be such that
$\N\cap\z(e)\subset \overline{{\cal O}(e)}$.
Then every nilpotent element of $\z(e;0)$ is contained in
$\z(e;0)\cap [e,\g(-2;\lambda_e)]$.
\end{lemma}
\pf Given a subvariety $X$ in $\g\cong {\mathbb A}^n$ and $x\in X$
we denote by $T_x(X)$ the tangent space to $X$ at $x$. We view
$T_x(X)$ as a linear subspace of $\g$. Our assumption on $G$
implies that the orbit map $G\rightarrow{\cal O}(e)$,
$\,g\mapsto(\text{Ad}\,g)\,e$ is separable; see [\cite{SpSt},
Chap.~I, Sect.~5]. Therefore, $T_{e}({\cal
O}(e))=T_e(\overline{{\cal O}(e)})=[\g,e].$ Since $\N\cap
\z(e)\subset \overline{{\cal O}(e)}$ we then have
$T_e(\N\cap\z(e))\subseteq [\g,e].$ By our earlier remarks, $$
T_e(\N\cap\z(e))=T_e\big(\N(\z(e;0)+\textstyle{\bigoplus_{i>0}}\,\z(e;i)\big)=
T_0(\N(\z(e;0))\oplus
T_e\big(\textstyle{\bigoplus_{i>0}}\,\z(e;i)\big). $$ This yields
$T_0(\N(\z(e;0))\subset [\g,e]$.

Now let $x$ be an arbitrary element in $\N(\z(e;0))$, a conical variety,
and consider the natural injection $\iota\colon Kx\hookrightarrow \N(\z(e;0))$.
Since $\iota$ is a closed immersion with $({\mathrm d}\iota)_0 (x)=x$ we get
$x\in T_0(\N(\z(e;0))$. But then $x\in [\g,e]\cap \z(e;0)=\z(e)\cap [e,\g(-2;\lambda_e)]$ as
claimed.
\qed
\smallskip

\noindent {\bf 3.2.} We say that a nilpotent element $e\in\g$ is
{\it almost distinguished} if the connected component of $Z_G(e)$
is a solvable group. For $p=0,$ this is consistent with our
definition in (2.1).
\begin{prop}
Let $G$ be as above and let $
e\in\N$ be such that  $\N\cap\z(e)\subset \overline{{\cal O}(e)}$.
Then $e$ is almost distinguished in $\g$.
\end{prop}
\pf Our arguments
will rely on various case-by-case considerations.
Let $V$ be a finite dimensional vector space over $k$.

\smallskip

\noindent (1) Suppose $G\,=\,\GL(V)$. This case is particularly
interesting for us since we have in mind an application to the
punctual Hilbert scheme ${\cal H}_n$ over $k$. Let
$V=V_1\oplus\cdots \oplus V_s$ be a decomposition of $V$ into a
direct sum of indecomposable $ke$-modules such that $\dim
V_1\ge\cdots\ge \dim V_s$. The partition of $\dim V$ associated
with $e$ is nothing but $(\dim V_1,\dim V_2,\ldots,\dim V_s)$.

Let $I\,=\,\big\{i\in\{1,\ldots,s-1\}\,|\,\dim V_i=\dim V_{i+1}
\big\}.$ Suppose $I\ne\emptyset$ and let $j$ be the smallest
integer in $I$. Let $d=\dim V_j=\dim V_{j+1}$ and  $M=V_j\oplus
V_{j+1}$. Let $v_1,\ldots, v_d$ and $w_1,\ldots,w_d$ be bases of
$V_j$ and $V_{j+1}$ such that $e(v_d)=e(w_d)=0$ and
$e(v_i)=v_{i+1}$, $\,e(w_i)=w_{i+1}$ for $i<d$. There is
$z\in{\frak gl}(M)$ such that $z(w_d)=0$, $\,z(v_i)=w_i$ for $1\le
i\le d$ and $\,z(w_i)=v_{i+1}$ for $1\le i<d$. For $i>d$,  put
$v_i=w_i=0$. Then $ze(v_i)=w_{i+1}=ez(v_i)$ and
$ze(w_i)=v_{i+2}=ez(w_i)$ for all $i$. In other words, the
endomorphisms $e_{\vert M}$ and $z$ commute. Note that $z$ acts on
$M$ as a nilpotent Jordan block of order $2d$. There is
$\hat{z}\in{\frak gl}(V)$ such that  $\hat{z}_{\vert M}=z$ and
$\hat{z}_{\vert V_i}=e_{\vert V_i}$ for all $i\not\in\{j,j+1\}$.
By construction, $\hat{z}\in\N\cap\z(e)$ and the partition of
$\dim V$ associated with $\hat{z}$ is not dominated by that of
$e$. Applying [\cite{Ge}] we now deduce that
$\hat{z}\not\in\overline{\GL(V)\,e}$. Since this contradicts our
assumption on $e$ it must be that $I=\emptyset$.

As a result, all parts of the partition of $\dim V$ associated with $e$ are distinct.
But then $e$ is almost distinguished in $\g$, by [\cite{SpSt}, Chap.~IV, Cor.~1.8(i)].
So the proposition holds in the present case.

\smallskip

\noindent
(2) We now consider the case where $G$ is  a  group of type $\mathrm{B}$,
$\mathrm{C}$ or $\mathrm{D}$.
Since $p$ is good for $G$ we have $p\ne 2$. Assume that $\dim V\ge 2$
and let $\Psi$ be a nondegenerate bilinear form on $V$ such that $\Psi(u,v)=
(-1)^{\varkappa}\,
\Psi(v,u)$ for all $u,v\in V$, where $\varkappa\in\{0,1\}$. Let $G(\Psi)$ be the
closed subgroup of $\SL(V)$ consisting of all
$g\in\SL(V)$ with $\Psi(g(u),g(v))=\Psi(u,v)$ for all $u,v\in V$. Then
$$\g(\Psi)\,:=\,\{x\in\sl(n)\,|\,\,\Psi(x(u),v)+\Psi(u,x(v))=0\ \, \text{for all }\, u,v\in V\}$$
is the Lie algebra of $G(\Psi)$. It is well-known that $G$ is isogenic to
$G(\Psi)$ for a suitable choice of $V$ and $\Psi$, and $\g\cong\g(\Psi)$ as
restricted Lie algebras. We thus may  identify $\g$ with $\g(\Psi)$ and view our
$e\in\N$ as a nilpotent endomorphism of $V$.

Let $(n_1\ge n_2\ge\cdots\ge n_s)$ be the partition of $\dim V$
associated with $e$. According to [\cite{SpSt}, Chap.~IV, (2.19)]
there is a direct sum decomposition $V=V_1\oplus\cdots\oplus V_s$
with $\dim V_i=n_i$ for all $i$, such that
\begin{itemize}
\item[1.] each $V_i$ is $e$-stable and
$e$ acts on $V_i$ as a nilpotent Jordan block of order $n_i$;
\smallskip
\item[2.] if $n_i+\varkappa$ is odd, then $\Psi$ is nondegenerate on $V_i$ and
$\Psi(V_i,V_k)=0$ for  $k\ne i$;
\smallskip
\item[3.] if $n_i+\varkappa$ is even, then $\Psi$ vanishes on
$V_i\times V_i$ and there exists a unique
$i^*=i\pm 1$ such that $n_{i^*}=n_{i}$, $\,\Psi$ is nondegenerate on
$V_i\oplus V_{i^*},$ and $V_k$ is orthogonal to $V_i\oplus V_{i^*}$
for  $k\not\in\{i,i^*\}$.
\end{itemize}
For $i\le \dim V$ we denote by $r(i)$ the number of $k$ with
$n_k=i$. Let $I_1$ (respectively, $I_2$) be the set of all
$i\in\{1,\ldots,s\}$ such that $n_i+\varkappa$ is odd
(respectively, even) and $r(n_i)\ge 3$ (respectively, $r(n_i)\ge
2$). Put $I=I_1\cup I_2$ and suppose $I\ne\emptyset$. Let $j$ be
the smallest integer in $I$ and set $d:=n_j=\dim V_j$.

\smallskip

\noindent (a) First suppose that $d+\varkappa$ is odd. Then
$r(d)\ge 3$ implying $\dim V_j=\dim V_{j+1}=\dim V_{j+2}=d$. As
$d+\varkappa$ is odd, the subspaces $V_j$, $V_{j+1}$ and $V_{j+2}$
are orthogonal to each other relative to $\Psi$. Let $M=V_i\oplus
V_{j+1}\oplus V_{j+2}.$ For $t\in\{0,1,2\}$, there is a basis
$m_{1,t},\ldots,m_{d,t}$ of $V_{j+t}$ such that $e(m_{d,t})=0$,
$\,e(m_{i,t})=m_{i+1,t}$ for $1\le i\le d-1$,  and
$$\Psi(m_{i,t},m_{k,t})=(-1)^{i-1}\,\delta_{i,\,d+1-k}\qquad\ \
(1\le i,k\le d);$$ see [\cite{SpSt}, Chap.~IV, (2.19)]. For $1\le
i\le d$, define
$u_i:=\frac{1}{\sqrt{2}}(m_{i,0}+\sqrt{-1}\,m_{i,2})$,
$\,v_i:=m_{i,1}$,
$w_i:=-\frac{1}{\sqrt{2}}(m_{i,0}-\sqrt{-1}\,m_{i,2})$ and let
$u_i=v_i=w_i=0$ for all $i>d$. Let $U$ and $W$ be the linear spans
of the $u_i$'s and $w_i$'s, respectively. By construction, these
subspaces are totally isotropic with respect to $\Psi$.
Furthermore, $M=U\oplus V_{j+1}\oplus W$ and $\Psi$ vanishes on
$(U+W)\times V_{j+1}$. Finally,
$\Psi(u_i,w_k)=-(-1)^{i-1}\,\delta_{i,\,d+1-k}=-\Psi(v_i,v_k)$ for
all  $i,k\le d$. There is $z\in{\frak gl}(M)$ such that
$z(u_i)=v_i$, $\,z(v_i)=w_i$ and $z(w_i)=u_{i+1}$ for all $i$. It
acts on $M$ as a nilpotent Jordan block of order $3d$. Similar to
part~1 of this proof one checks that the endomorphisms $e_{\vert
M}$ and $z$ commute. Next observe that
$\Psi(z(u_i),v_k)=\Psi(v_i,v_k)=-\Psi(u_i,w_k)=-\Psi(u_i,z(v_k))$
and
\begin{eqnarray*}
\Psi(z(w_i),w_k)&=&\Psi(u_{i+1},w_k)=-(-1)^{i}\,\delta_{i+1,d+1-k}=
(-1)^{d+\varkappa-i}\,\delta_{i+1,d+1-k}\\
&=&(-1)^{k+\varkappa}\,\delta_{k+1,d+1-i}=-(-1)^{\varkappa}\,\Psi(u_{k+1}, w_{i})
=-\Psi(w_i,z(w_k)).
\end{eqnarray*}
Combined with our earlier remarks this shows that the endomorphism
$z$ is skew-adjoint with respect to $\Psi_{\vert M\times M}$.
Property~2 of our direct sum decomposition now ensures that there
exists $\hat{z}\in\g=\g(\Psi)$ such that $\hat{z}_{\vert M}=z$ and
$\hat{z}_{\vert V_i}=e_{\vert V_i}$ for all
$i\not\in\{j,j+1,j+2\}$. By construction, $\hat{z}\in\N\cap\z(e)$
and the partition of $\dim V$ associated with $\hat{z}$ is not
dominated by that of $e$. So $\hat{z}\not\in\overline{\GL(V)\,e}$,
by [\cite{Ge}]. But then $\hat{z}\not\in\overline{G\,e}$, a
contradiction.

\smallskip

\noindent (b) Now suppose  $d+\varkappa$ is even. This case is
harder as we do not have much room for manoeuvre here (indeed, it
may happen that $r(d)=2$). Set $M:=V_j\oplus V_{j^*}=V_j\oplus
V_{j+1}$.

By property~3 of our decomposition, the subspaces $V_j$ and
$V_{j+1}$ are totally isotropic with respect to $\Psi$. According
to [\cite{SpSt}, Chap.~IV, (2.19)] there exist bases
$v_1,\ldots,v_d$ and $v'_1,\ldots,v'_d$ of $V_j$ and $V_{j+1}$
respectively, such that $e(v_d)=e(v'_d)=0$, $\,e(v_{i})=v_{i+1},$
$\,e(v'_i)=v'_{i+1}$ for $1\le i<d$, and
$$\Psi(v_i,v'_k)=(-1)^{i-1}\,\delta_{i,\,d+1-k} \qquad\, (1\le
i,k\le d).$$ Since $d+\varkappa$ is even we have
$\Psi(v'_i,v_k)=-\Psi(v_i,v'_k).$ For $i>d$, put $v_i=v'_i=0$.

Let ${\cal A}$ be the algebra $k[X]/(X^d)$ and $x$ be the image of $X$ in $\cal A$.
Denote by $C$ the centraliser of $e_{\vert M}$ in $\text{End}(M)$.
Since $e_{\vert M}$ acts on $M$ as a direct sum of two Jordan blocks of order $d$,
there exists an isomorphism of associative algebras $\varphi\colon\,C\,\stackrel{\sim}
{\longrightarrow}\,\Mat_2({\cal A})$ such that $\varphi(e_{\vert M})\,=\,
\text{diag}\,(x,x)$. Specifically, if
$\phi_{hk}(x)$ is the $(h,k)$th entry of $\varphi(c)$ for $c\in C$,
then $
c(v_1)\,=\,\phi_{11}(e)\,v_1+\phi_{21}(e)\,v'_2$ and
$c(v'_1)\,=\,\phi_{12}(e)\,v_1+\phi_{22}(e)\,v'_2.$

Now let $z\in C$ be such that
$\varphi(z)=\left({x\atop {\,x^2}}\, {1\,\atop {\,x\,}}\right)$.
We claim that $z$ is skew-adjoint with respect
to $\Psi_{\vert M\times M}$. To prove this we first observe that
$z(v_i)=v_{i+1}+v'_{i+2}$ and $z(v'_i)=v'_{i+1}+v_i$ for all $i$. Then
\begin{eqnarray*}
\Psi(z(v'_i),v'_{k})&=&\ \Psi(v_{i},v'_k)\ \, =\ -\Psi(v'_{i},v_k)\ \
=-\Psi(v'_i,z(v'_k)),\\
\Psi(z(v_i),v_{k})&=&\Psi(v'_{i+2},v_k)=-\Psi(v_{i+2},v'_k)=-\Psi(v_i,z(v_k)),\\
\Psi(z(v_i),v'_{k})&=&\Psi(v_{i+1},v'_k)=-\Psi(v'_{i},v_{k+1})
=-\Psi(v_i,z(v'_k)),
\end{eqnarray*}
hence the claim. Next observe that
$$\varphi(z)^2=\left(\begin{array}{cc} x& 1\\
x^2&x\end{array}\right)\cdot\left(\begin{array}{cc} x& 1\\
x^2&x\end{array}\right)\,=\,\left(\begin{array}{cc} 2x^2& 2x\\
2x^3&2x^2\end{array}\right)=2x\,\varphi(z).$$
Easy induction on $r$ now shows that
$\varphi(z)^{r}=(2x)^{r-1}\varphi(z)$ for all $r\in\mathbb N$.
But then $$\phi(z)^d\,=\,(2x)^{d-1}\left(\begin{array}{cc} x& 1\\
x^2&x\end{array}\right)=\,2^{d-1}\left(\begin{array}{cc} 0& x^{d-1}\\
0&0\end{array}\right)\ne\,0,\qquad \varphi(z)^{d+1}=0,$$
which implies that $z^{d}\ne 0$ and $z^{d+1}=0$.
So $z$ is nilpotent and at least one part of the partition of $\dim M$ associated
with $z$ equals $d+1$. Property~3
of our decomposition ensures that there exists
$\hat{z}\in\g$ with
$\hat{z}_{\vert M}=z$ and $\hat{z}_{\vert V_i}=e_{\vert V_i}$ for all $i\not\in\{j,j+1\}$.
It is clear from the discussion above that $\hat{z}\in\N\cap\z(e)$ and the partition
of $\dim V$ associated with $\hat{z}$ is not dominated by that of $e$. Arguing as
at the end of part~2a of this proof  we now deduce that
$\hat{z}\not\in\overline{G\,e}$, a contradiction.

When combined, parts~2a and 2b show that $I=\emptyset$. In other words,
$n_i+\varkappa$ is odd and $r(n_i)\le 2$ for all $i\le s$. But then
[\cite{SpSt}, Chap.~IV, (2.25)] yields that
the connected centraliser $Z_G(e)^\circ$ is solvable (one should keep
in mind that the group $\SO(2)$ is a torus). Thus $e$ is almost distinguished
when $G$ is of type $\mathrm B$, $\mathrm C$ or $\mathrm D$.

\smallskip

\noindent
(3) Now suppose $G$ is a group of type $\mathrm E$.
We are going to rely on some results proved in [\cite{P1}]. Adopt
Dynkin's labelling of nilpotent orbits; see [\cite{Ca}, Chap.~13].
Let $G_{\mathbb C}$ be a simple algebraic group over $\mathbb C$
of the same type as $G$ and let
$e_{\mathbb C}$ be a nilpotent element in
$\Lie\, G_{\mathbb C}$ whose $G_{\mathbb C}$-orbit
has the same labelling as ${\cal O}(e)\subset\g$.
It follows from the main result in [\cite{Spa}] (and also from
[\cite{P02}, Sect.~2]) that $$\dim\z(e)=\dim_{\mathbb C}\,\z(e_{\mathbb C}),\quad\ \,
\dim \z(e;0)=
\dim_{\mathbb C}\, (\z(e_{\mathbb C})/\text{nil}\,\z(e_{\mathbb C})).$$
So we can use Elashvili's tables [\cite{Ca}, pp.~401--407]
for computing $\dim \z(e)$ and $\dim \z(e;0)$.

\smallskip

\noindent
(a) Suppose $e$ is not almost distinguished in $\g$
(then $e$ is not distinguished in $\g$ either).
Since the connected component of
$Z_G(e)/R_u(Z_G(e))\cong Z(\lambda_e)\cap Z_G(e)$
is not a torus, by our assumption on $e$, and since
 $\z(e;0)=\Lie\big(Z(\lambda_e)\cap Z_G(e)\big)$,
by [\cite{Jan}, Sect.~5] and [\cite{P02}, Sect.~2], it must be
that $\N(\z(e;0))\ne\{0\}$. Then $\z(e;0)\cap
[e,\g(-2;\lambda_e)]\ne\{0\}$, by Lemma~3.1. So the linear map
$(\text{ad}\,e)^2\colon\,\g(-2;\lambda_e)\longrightarrow\g(2;\lambda_e)$
is not bijective. The computations in [\cite{P1}] now show that
$p\in\{5,7\}$ and $e$ is regular  in $\text{Lie}\,L$ where $L$ is
a Levi subgroup of $G$ with a factor of type $A_{p-1}$. There are
eight such cases in all, and in each of them we have $e^{[p]}=0$.

If $G$ is of type ${\mathrm E}_6$ (respectively, ${\mathrm E}_7$)
and ${\cal O}(e)$ is labelled by
$\mathrm A_4\times A_1$, then $\z(e;0)$ is $1$-dimensional
(respectively, $2$-dimesional); see [\cite{Ca}, pp.~402, 404].
Combined with our earlier remarks this implies that
the connected component of
$Z(\lambda_e)\cap Z_G(e)$ is a torus. Hence
${\cal O}(e)$ is not of the above type.
If ${\cal O}(e)$
is labelled by ${\mathrm A}_{p-1}$ and $G$ is not of type
$\mathrm E_7$ when $p=7$, then $\z(e)$
meets the orbit
labelled by ${\mathrm A}_{p-1}\times \mathrm A_1$, call it ${\cal O}_1$.
Elashvili's tables along with
our earlier remarks assure that $\dim {\cal O}_1>\dim {\cal O}(e)$. Then
$\N\cap \z(e)\not\subset\overline{{\cal O}(e)}$, hence ${\cal O}(e)$ is not of that type.

These observations settle the case where $G$ is of type $\mathrm E_6$ and
leave us with just three orbits, one in characteristic $5$
and two in characteristic $7$. Unfortunately,
the conclusion of Lemma~3.1 does hold for  these orbits
and no further reduction is readily available.
So we have to work harder here, and
our plan  will be
to exhibit, in each of the remaining cases, an element $z\in\N\cap\z(e)$ with
$z^{[p]}\ne 0$. Since $e^{[p]}=0$ and the $[p]$th power map
on $\g$ is a $G$-equivariant morphism, this will imply that
$\N\cap\z(e)\not\subset\overline{{\cal O}(e)}$.
We adopt Bourbaki's numbering of simple roots [\cite{Bou}, Tables~I--IX]
and the notation of [\cite{P1}]. The group scheme $(SL_2)_H$ from
[\cite{P1}, (2.26)] will be denoted by $\cal G$.

\smallskip

\noindent
(b) Suppose $p=7$ and $G$ is of type $\mathrm E_8$. Then ${\cal O}(e)$ is labelled
by $\mathrm A_6\times A_1$. So we may assume that $L=L_J$, where
$J=\{1,2,4,5,6,7,8\}$, and $e=\sum_{i\in J}e_{\alpha_i}$ (note that $I=J$ in this case).
From [\cite{P1}, (2.27)] we know that the projective $\cal G$-module $P_{0,0}$
is a direct summand of both $\g_J(-2)$ and $\g_J(2)$. So
$\z(e;0)\cap\g_J(\pm2)\ne \{0\}$. On the other hand, $\dim \z(e;0)=3$, by
[\cite{Ca}, p.~406]. Then $\z(e;0)$ is a $3$-dimensional reductive Lie algebra with
$\N(\z(e;0))\ne \{0\}$. Therefore, $\z(e;0)\cong\sl(2)$ as restricted Lie algebras.

Obviously, $Z(\lambda_e)\cap Z_G(e)$ contains $Z(L_J)$, the centre
of $L_J$. Combined with the preceding remark this shows that there
exist nonzero $e_0\in\g_J(2)\cap\g(0;\lambda_e)$,
$f_0\in\g_J(-2)\cap\g(0;\lambda_e)$ and $h_0\in\Lie\,Z(L_J)$ such
that $(e_0,h_0,f_0)$ is an $\sl_2$-triple in $\z(e;0)$. Since each
$\g_J(\ell)$ with $\ell\ne 0$ is an irreducible $L_J$-module, by
[\cite{P1}, (2.27)], the endomorphism $\text{ad}\,h_0$ acts on
$\g_J(\ell)$ as $\ell\,\text{id}$.

By [\cite{P1}, (2.27)], we have that $\g_J(\pm 4)\cong V_{p-1,0}$
as $\cal G$-modules. In particular, the subspaces
$\z(e)\cap\g_J(\pm 4)$ are $1$-dimensional. Fix a nonzero
$a\in\z(e)\cap\g_J(4)$, a multiple of $e_{\tilde{\alpha}}$, and
put $z=f_0+a$. We claim that $z\in\N$ and $z^{[7]}\ne 0$.

For $m_1,\ldots, m_{s}\in \mathbb Z_+$ and $n_1,\ldots,n_s
\in\mathbb N$,  set $$[a^{m_1}f_0^{n_1}\cdots\,
a^{m_s}f_0^{n_s}a]:=
(\text{ad}\,a)^{m_1}(\text{ad}\,f_0)^{n_1}\cdots\,(\text{ad}\,a)^{m_s}
(\text{ad}\,f_0)^{n_s}(a).$$ To prove the claim we first observe
that $a\in\g(6;\lambda_e)$ is a primitive vector of weight $4$ for
$\z(e;0)$. Hence $[f_0^4a]\ne 0$, by a standard $\sl_2$-argument.
Since $[f_0^4a]$ belongs to $\g_J(-4)\cap\z(e)$, a $1$-dimensional
subspace, it must be that
\begin{eqnarray}
[f_0^4a]\in k^*e_{-\delta}\quad\text{ where } \ \,
\textstyle{\delta={2\atop{}}{4\atop{}}{5\atop
2}{4\atop{}}{3\atop{}}{2\atop{}} {1\atop{}}}.
\end{eqnarray}

By [\cite{Jac}, Chap.~V, Sect.~7],  $z^{[7]}=a^{[7]}+f_0^{[7]}+
\sum_{i=1}^{6} s_i(a,f_0)$
where $$\text{ad}(ta+f_0)^6(a)=\textstyle{\sum_{i=1}^6}
\,is_i(a,f_0)t^{i-1}\quad\ \ \, (\forall\,t
\in k).$$ Since $\g_J(6)$ is zero, we have $[af_0\,a]=[a^2f_0^3a]=0$. Then
$[af_0^2 a]=0$, by the Jacobi identity. It follows that
\begin{eqnarray}
\text{ad}(ta+f_0)^6(a)=\,\text{ad}(ta+f_0)^3([f_0^3a])=\,\text{ad}(ta+f_0)^2(t[af_0^3a]
+[f_0^4a]).
\end{eqnarray}
As $\g_J(-6)$ is zero, $[f_0^5a]=0$ necessarily holds. Since both
$[a^2f_0^4a]$ and $[af_0af_0^3a]$ belong to
$\g_J(4)\cap\g(18;\lambda_e)$, we also have
$[a^2f_0^4a]=[af_0af_0^3a]=0.$ So (6) yields
\begin{eqnarray}
\text{ad}(ta+f_0)^6(a)=\,t[ta+f_0,\,[af_0^4a]+[f_0af_0^3a]]=\,t([f_0af_0^4a]+
[f_0^2af_0^3a]).
\end{eqnarray}
Now let $F=\text{ad}\,f_0$ and $A=\text{ad}\,a$. Since
$[f_0^5a]=0$ we have $(\text{ad}\,F)^5(A)$ forcing
$$(F^5A-5F^4AF+10F^3AF^2-10F^2AF^3 +5FAF^4-AF^5)(a)=0.$$ Together
with our observations above this gives
$[f_0a^4f_0\,a]=2[f_0^2af_0^3a].$ By (7),
$$\text{ad}(ta+f_0)^6(a)=\,
=\,3t[f_0^2af_0^3a]\,=\,-4t[f_0^2af_0^3a]\,=\, 2s_2(a,f_0)t.$$ So
$s_i(a,f_0)=0$ for $i\ne 2$.  Also, $f_0^{[7]}=a^{[7]}=0,$ because
$\z(e;0)\cong\sl(2)$ as restricted Lie algebras. It follows that
\begin{eqnarray}
z^{[7]}=\,s_2(a,f_0)\,=-2[f_0^2af_0^3a]\in\,\g_J(-2)\subset\N.\end{eqnarray}
According to [\cite{Bou}, Table~VIII],
$\tilde{\alpha}-\delta$ is a root. So (5) yields $[af_0^4a]\ne 0$.
This, in turn, implies that
$[af_0^3a]\ne 0$. Indeed,
otherwise $$[af_0^4a]=[[af_0], [f_0^3a]]=[f_0,[[a,f_0],[f_0^2a]]]
=[f_0,[af_0^3a]]=0,$$ a contradiction (here we used the Jacobi identity and
the fact that $[af_0^2a]=0$).  Since
$[e_0,[af_0^3a]]\in k[af_0^2a]$, we obtain that
$[af_0^3a]$ is
a primitive vector of weight $2$ for $\z(e;0)$. But then $[f_0^2af_0^3a]\ne 0$,
by a standard $\sl_2$-argument. Now(8) shows that $z\in\N$ and
$z^{[7]}\ne 0$, as claimed. This settles the present case.

\smallskip

\noindent (c) Retain the assumptions of part~3b. Since $G$ is of
type $\mathrm E_8$, it is both adjoint and simply connected. In
particular, there is a unique involution $\tau\in T\subset G$ with
the property that that $(\text{Ad}\,\tau)\,x=(-1)^\ell x$ for all
$x\in\g_J(\ell)$ and all $\ell\in\mathbb Z$. Let
$G^\tau=Z_G(\tau)$ and let $\g^\tau$ be the fixed point algebra of
$\text{Ad}\,\tau$. By a result of Steinberg, $G^\tau$ is a
connected reductive group; see [\cite{SpSt}, Chap.~II, \S~4]. By
[\cite{Bo}, Chap.~3, Sect.~9], we have $\g^\tau=\Lie\,G^\tau$.

Clearly, $T\subset G^\tau$. Straightforward calculation shows that
$\dim\,\g^\tau=136$. Since all maximal root subsystems in the root
system of type $\mathrm E_8$ are well-known (see [\cite{Bou},
Chap.~VI, \S~4, Exercise~4] for example), one finds out without
difficulty that $G^\tau$ is a semisimple group of type $\mathrm
E_7\times A_1$. Let $G_1$ and $G_2$ be the simple components of
$G^\tau$ of type $\mathrm E_7$ and $\mathrm A_1$, respectively,
and $\g_i=\Lie\,G_i$. The Lie algebras $\g_1$ and $\g_2$ are
simple, hence $\g^\tau=\,\g_1\oplus \g_2$, a direct sum of
restricted Lie algebras.

Given $x\in\g^\tau$ we let $x_i$ be the component of $x$ in $\g_i$.
Then $x^{[7]}=x_1^{[7]}+x_2^{[7]}$.
Since $\g_2\cong\sl(2)$, we also have that $x_2^{[7]}=0$ for any $x\in\g^\tau\cap\N$.
Note that $e,z\in\g^\tau$.
The preceding remark shows that $e_1^{[7]}=0$ and
$z_1^{[7]}\ne 0$. Obviously,  $[e_1,z_1]=0$.

Next observe that $P_J\cap G^\tau$ is a parabolic subgroup of
$G^\tau$ and $L_J$ is a Levi subgroup in it. From this it follows
that $L_J\cap G_1$ is a Levi subgroup of type $A_6$ in $G_1$.
Since $e$ is regular in ${\frak l}_J$, its component $e_1$ must be
regular in ${\frak l}_J\cap\g_1\,=\,\Lie\,(L_J\cap G_1)$. By the
discussion above, $z_1\not\in\overline{(\text{Ad}\,G_1)\,e_1}$.
This settles our second remaining case in characteristic $7$.

\smallskip

\noindent (d) Now suppose $p=5$. Then $G$  is of type $\mathrm
E_7$ and ${\cal O}(e)$ is labelled by $\mathrm A_4\times A_2$; see
part~3a. So we may assume that $L=L_J$, where $J=\{1,2,3,4,6,7\}$,
and $e=\sum_{i\in J}e_{\alpha_i}$. Since
$\tilde{\alpha}={2\atop{}}{3\atop{}}{4\atop
2}{3\atop{}}{2\atop{}}{1\atop{}}$, we have
$M_{J,+}=\,\g_J(1)\oplus \g_J(2)\oplus\g_J(3)$. Using [\cite{Bou},
Table~VI] we observe that $\g_J(1)$, $\g_J(2)$ and $\g_J(3)$ are
irreducible modules over $(L_J,L_J)\cong \SL_5(k)\times\SL_3(k)$
with highest weights $\omega_3^J+\omega_7^J$,
$\omega_2^J+\omega_6^J$ and $\omega_1^J$, respectively. As a
consequence, $$\g_{J}(1)\cong \,(\stackrel{2}{\wedge} N_5)\otimes
N_3,\quad \ \ \ \g_J(2)\cong\, N_5^*\otimes N_3,\quad \ \ \
\g_J(3)\cong\,N_5.$$ As in [\cite{P1}, (2.27)] this yields that
each $\g_J(\ell)$ with $\ell\ne 0$ is a projective $\cal
G$-module. Direct calculation based on [\cite{Bou}, Table~VI] and
[\cite{P1}, Lemma~2.9]  gives $m_1(e)=8$, $\,m_2(e)=6$,
$\,m_3(e)=4$ showing that $2p-2$ is the largest $\lambda_e$-weight
of $M_{J,+}$. Then [\cite{P1}, (2.26)] says that all
indecomposable summands of the $\cal G$-modules $M_{J,\pm}$ are of
the form $V_{p-1, 0}$ or $P_{m, 0}$ for $m\le p-2$. Moreover,
$P_{0,0}$ occurs in $\g_J(\pm 1)$.

On the other hand, $\dim \z(e;0)=3$, by [\cite{Ca}, p.~406].
Arguing as in part~3b we now obtain that there are nonzero
$e_0\in\g_J(1)\cap\g(0;\lambda_e)$,
$f_0\in\g_J(-1)\cap\g(0;\lambda_e)$ and $h_0\in\Lie\,Z(L_J)$ such
that $(e_0,h_0,f_0)$ is an $\sl_2$-triple in $\z(e;0)$. Since each
$\g_J(\ell)$ with $\ell\ne 0$ is $L_J$-irreducible,
$\text{ad}\,h_0$ acts on $\g_J(\ell)$ as $2\ell\,\text{Id}$
(because $e_0\in \g_J(1)$ and $[h_0,e_0]=2e_0$). As in part~3b,
$\z(e;0)\cong\,\sl(2)$ as restricted Lie algebras.

Since $(\g_J(-3))^*\cong\g_J(3)\cong N_5$ as $(L_J,L_J)$-modules
and $e$ is regular nilpotent in ${\frak l}_J$, it must be that
$\g_J(\pm 3)\cong V_{p-1,0}$ as $\cal G$-modules. In particular,
the subspaces $\z(e)\cap\g_J(\pm 3)$ are $1$-dimensional. Fix a
nonzero $a\in\z(e)\cap\g_J(3)$, an element in $\g(4;\lambda_e)$,
and put $z=f_0+a$. We claim that $z\in\N$ and $z^{[5]}\ne 0$.

By [\cite{Jac}, Chap.~V, Sect.~7], $z^{[5]}=a^{[5]}+f_0^{[5]}+
\sum_{i=1}^{4} s_i(a,f_0)$ where
$$\text{ad}(ta+f_0)^4(a)=\textstyle{\sum_{i=1}^4}
\,is_i(a,f_0)t^{i-1}\quad\ \ \, (\forall\,t \in k).$$ Since
$\g_J(5)$ is zero, we have $[af_0\,a]=0$. Then $[af_0^2 a]=0$, by
the Jacobi identity. Furthermore, $[af_0^3a]=0$ as $\g_J(3)\cap
\g(8;\lambda_e)$ is zero.  Hence
$$\text{ad}(ta+f_0)^4(a)=\,\text{ad}(ta+f_0)([f_0^3a])=\,[f_0^4a]=s_1(a,f_0).$$
This means that $s_i(a,f_0)=0$ for $i>1$. Since
$a^{[5]}=f_0^{[5]}=0$, we now get
\begin{eqnarray}
z^{[5]}=[f_0^4a]\in\,\g_J(-1)\subset\N.
\end{eqnarray}
Since $a$ is a primitive vector of weight $1$ for $\z(e;0)$, the element
$(\text{ad}\,e_0)^2([f_0^4a])$ is a nonzero multiple of $[f_0^2a]$.
Due to (9) we are thus reduced to show that $[f_0^2a]\ne 0$.

Since $\beta_J^5={1\atop{}}{2\atop{}}{2\atop 1}{1\atop{}}{1\atop{}}{1\atop{}}$ and
$b_J^5=10$, it follows from [\cite{P1}, Lemma~2.9] and [\cite{Bou}, Table~VI]
that the subspace $\g_J(-1)\cap\g(0;\lambda_e)$ is spanned by
$e_{-\gamma_1},\ldots,e_{-\gamma_6}$ where
\begin{eqnarray*}
\gamma_1&=&\textstyle{{0\atop{}}{0\atop{}}{1\atop 1}{1\atop{}}{1\atop{}}{1\atop{}}},
\qquad
\gamma_2\,=\,\textstyle{{0\atop{}}{1\atop{}}{1\atop 0}{1\atop{}}{1\atop{}}{1\atop{}}},
\qquad
\gamma_3\,=\,\textstyle{{1\atop{}}{1\atop{}}{1\atop 0}{1\atop{}}{1\atop{}}{0\atop{}}},
\\
\gamma_4&=&\textstyle{{0\atop{}}{1\atop{}}{1\atop 1}{1\atop{}}{1\atop{}}{0\atop{}}},
\qquad
\gamma_5\,=\,\textstyle{{1\atop{}}{1\atop{}}{1\atop 1}{1\atop{}}{0\atop{}}{0\atop{}}},
\qquad
\gamma_6\,=\,\textstyle{{0\atop{}}{1\atop{}}{2\atop 1}{1\atop{}}{0\atop{}}{0\atop{}}}.
\end{eqnarray*}
Let $\beta=\sum_{i=1}^7\alpha_i$, a positive root. We may assume without
loss of generality that
\begin{eqnarray*}
e_{-\gamma_1}&=&[e_{\alpha_3},[e_{\alpha_1},e_{-\beta}]],\quad
e_{-\gamma_2}\,=\,[e_{\alpha_2},[e_{\alpha_1},e_{-\beta}]],\quad
e_{-\gamma_3}\,=\,[e_{\alpha_2},[e_{\alpha_7},e_{-\beta}]]\\
e_{-\gamma_4}&=&[e_{\alpha_7},[e_{\alpha_1},e_{-\beta}]],\quad
e_{-\gamma_5}\,=\,[e_{\alpha_6},[e_{\alpha_7},e_{-\beta}]],\quad
[e_{\alpha_4},e_{-\gamma_6}]\,=\,[e_{\alpha_1},e_{-\gamma_5}].
\end{eqnarray*}
Let $f_0=\sum_{i=1}^6 s_i\,e_{-\gamma_i}$ where $s_i\in k$. Recall
that $[e,f_0]=0$ and $e=\sum_{i\ne 5}e_{\alpha_i}$. Since
$[e_{\alpha_{2}},e_{\alpha_{3}}]=[e_{\alpha_3}, e_{\alpha_7}]=
[e_{\alpha_{2}},e_{\alpha_6}] =0$, this yields $
s_1+s_2\,=\,s_1+s_4\,= \,s_3+s_5\,=\,0.$ Since $e_{\alpha_1},\,
e_{\alpha_2},\, e_{\alpha_7}$ pairwise commute we also have
$s_2+s_3+s_4=0.$ Since $[e_{\alpha_1},e_{\alpha_6}]=0$, our choice
of $e_{-\gamma_6}$ yields $s_4+s_5+s_6=0$.

Together these relations show that the $s_i$'s are  all nonzero.
Since $\g_J(3)\cong V_{p-1,0}$ as $\cal G$-modules, we have $a\in
\,k^*e_{\tilde{\alpha}}$. In view of [\cite{Bou}, Table~VI] we
have $$[f_0,e_{\tilde{\alpha}}]\in \{\pm
s_3\,e_{{1\atop{}}{2\atop{}}{3\atop 2}{2\atop{}}{1\atop{}}
{1\atop{}}}\pm \,s_5\,e_{{1\atop{}}{2\atop{}}{3\atop
1}{2\atop{}}{2\atop{}} {1\atop{}}}\}. $$ Using [\cite{Bou},
Table~VI] one now observes without difficulty that $[f_0^2a]$ is a
nonzero linear combination of $e_{{0\atop{}}{1\atop{}}{2\atop
1}{1\atop{}}{1\atop{}}{1\atop {}}},\,\,$
$e_{{1\atop{}}{1\atop{}}{1\atop 1}{1\atop{}}{1\atop{}}{1\atop
{}}},\,\,$ $e_{{1\atop{}}{1\atop{}}{2\atop
1}{1\atop{}}{1\atop{}}{0\atop {}}}.$ This proves the proposition
in the case where $G$ is of type $\mathrm E$.

\smallskip

\noindent (4) Next suppose $G$ is of type $\mathrm F_4$. Then
$G\,=\,\text{Aut}(\g)$. Let $\tilde{G}$ be a simply connected
group of type $\mathrm E_6$ and $\tilde{\g}\,=\,\Lie\,\tilde{G}$.
As in [\cite{P1}, (2.28)], we regard $\g$ as a subalgebra of
$\tilde{\g}$. More precisely, we assume that
$\g\,=\,\tilde{\g}^\sigma$ where $\sigma$ is the involution in
$\text{Aut}(\tilde{\g})$ swapping $e_{\pm \alpha_i}$ and $e_{\pm
\alpha_{-i+7}}$ for $i=1,2$ and fixing $e_{\pm\alpha_i}$ for
$i=2,4$. Note that $\sigma$ permutes the set of indices
$\{1,2,\ldots,6\}$.

Suppose $e$ is conjugate to a regular nilpotent element in  a
standard Levi subalgebra of $\g$. Then a $G$-conjugate of  $e$ is
regular in a standard Levi subalgebra ${\frak l}_{\tilde{J}}$ of
$\tilde{\g}$ such that $\tilde{J}^\sigma=\tilde{J}$; see
[\cite{P1}, (2.28)]. Hence it can be assumed in the present case
that $e=\sum_{i\in\tilde{J}}e_{\alpha_i}$. Let $\tilde{\lambda}$
denote the image of $\lambda_{\tilde{J},\tilde{J}}\in
X_*(\tilde{G})$ under the natural embedding
$X_*(\tilde{G})\hookrightarrow X_*(\text{Aut}(\tilde{\g}))$. It is
easy to see that $\sigma$ fixes $\tilde{\lambda}$.

The construction of optimal $1$-parameter subgroups
in [\cite{P02}, Sect.~2] shows that we can take as
$\lambda_{e}$ the image of $\tilde{\lambda}$ in $X_*(G)=
X_*(\text{Aut}(\tilde{\g}^\sigma))$.
Since ${\frak l}_{\tilde{J}}$ is $\sigma$-stable and $p>3$,
the root system of ${\frak l}_{\tilde{J}}$ has no components of type
${\mathrm A}_{p-1}$. Calculations in [\cite{P1}]  then yield that the map
$(\text{ad}\,e)^2\colon\,\tilde{\g}(-2;\tilde{\lambda})\,\longrightarrow\,
\tilde{\g}(2;\tilde{\lambda})$ is bijective. It follows that so is the map
$(\text{ad}\,e)^2\colon\,\g(-2;\lambda_e)\,\longrightarrow\,
\g(2;\lambda_e)$.

If $e$ is not of the above type, then it follows from
[\cite{P1}, (2.28)] and the proof of [\cite{P1}, Lemma~2.7] that
the map $(\text{ad}\,e)^2\colon\,\g(-2;\lambda_e)\,\longrightarrow\,
\g(2;\lambda_e)$ is bijective. Arguing as in part~3a we now obtain
that in any event $e$ is almost distinguished in $\g$.
\smallskip

\noindent (5) Finally, suppose $G$ is of type $\mathrm G_2$. We
may assume that $e$ is not distinguished in $\g$. Note that $e\ne
0$ and $e$ is not conjugate to a long root vector in $\g$
(otherwise $\z(e)$ would contain a regular nilpotent element in
$\g$). By [\cite{Ca}, p.~401], we are now left with the orbit
labelled by $\tilde{\mathrm A}_1$. So suppose $e=e_{\alpha_2}$.
Then $I=J=\{2\}$ and we may assume further that
$\lambda_e=\alpha_2^\vee$. It is easily seen that in this case
$\z(e;0)\cong\sl(2)$ and $\dim\,\g(\pm 2;\lambda_e)=1$. But then
Lemma~3.1 shows that this case is impossible, completing the proof
of the proposition. \qed

\smallskip

\noindent
{\bf 3.3.}
Proposition~3.2 enables us now to establish
a modular version of Proposition~2.1.
\begin{prop}
Under the present assumption on $G$, each irreducible component of
${\frak C}^{\mathrm nil}(\g)$ is of the form ${\frak C}(e)$ where
$e$ is almost distinguished in $\g$.
\end{prop}
\begin{pf}
By [\cite{HSp}],
the nilpotent orbits in reductive Lie algebras are finite in number
regardless of $p$.  Therefore,
our arguments  in parts~1 and 2 of the proof of Proposition~2.1 apply
to {\it any} reductive Lie algebra.
They show that each irreducible component of $\C(\g)$ is of the
form ${\frak C}(e)$ where $e\in\N$ is such that
$\N\cap\z(e)\subset \overline{{\cal O}(e)}$. According to Proposition~3.2,
each such $e$ is almost distinguished in $\g$.
\end{pf}

\noindent {\bf 3.4.} Recall that all our observations in (2.2) are
valid in the present setting. With the notation of (2.2), we pick
$(I,J)$ in ${\cal P}(\Pi)$ and let $e\in{\frak
l}_{I}(2;\lambda_{I,J})$ be such that the orbit
$(\text{Ad}\,P_{I,J})\,e$ is open in ${\frak u}_{I,J}$. By
Lemma~2.2, $P_J$ is a distinguished parabolic subgroup of $G$; so
$\lambda_J\in X_*(T)$ is well-defined. By [\cite{Ca},
Prop.~5.8.5],  ${\frak p}_J\,=\,\bigoplus_{i\ge
0}\,\g(2i;\lambda_J)$ and the subspace $\g(2;\lambda_J)$ contains
a Richardson element of ${\frak p}_J$. As in (2.3), we pick such
an element and call it $\tilde{e}$. It follows from the definition
of $\lambda_J$ that
\begin{eqnarray}
\g(i;\lambda_J)\cap
{\frak l}_{I}\,=\,{\frak l}_{I}(i;\lambda_{I,J})  \ \quad \qquad(\forall\,i\in{\mathbb Z}).
\end{eqnarray}
In particular, $e\in\g(2;\lambda_J)$ so that
$\text{Int}\,\lambda_{J}(k^*)$ acts on $Z_G(e)$. As a result,
$\text{Ad}\,\lambda_J(k^*)$ acts on $\Lie\, R_u(Z_G(e))$ and on
each weight space $\g(i;\lambda_{I,J})$ of $\lambda_{I,J}(k^*)$.
\begin{lemma}
The torus ${\mathrm Ad}\,\lambda_J(k^*)$ has no zero weight
on $\Lie\,R_u(Z_G(e))$.
\end{lemma}
\begin{pf}
It follows from [\cite{Jan}, Sect.~5] and [\cite{P02}, Sect.~2] that
$$\Lie\,R_u(Z_G(e))\,=\,\textstyle{\bigoplus_{i>0}}\ \g(i;\lambda_{I,J})\cap\z(e),$$
while from the definition of $\lambda_J$ it is immediate
that $\g(0;\lambda_J)={\frak l}_J\subseteq {\frak l}_I$. Then
\begin{eqnarray*}
\Lie\,R_u(Z_G(e))\cap\g(0;\lambda_J)&=&
\textstyle{\bigoplus_{i>0}}\ \g(i;\lambda_{I,J})\cap{\frak l}_J\cap\z(e)
\\&\subseteq &\textstyle{\bigoplus_{i>0}}\ {\frak l}_I(i;\lambda_{I,J})\cap{\frak l}_J
\,=\,\{0\}, \end{eqnarray*}
in view of (10) and the equality ${\frak l}_J=\,{\frak l}_I(0;\lambda_{I,J})$.
The result follows.
\end{pf}

\noindent {\bf 3.5.} Given a linear algebraic group $H$ we denote
by ${\cal U}(H)$ the unipotent variety of $H$ and put ${\cal
U}={\cal U}(G)$. Note that ${\cal U}\subset (G,G)$. Since the
$\sl_2$-theory has its limitations, we have to modify our
constructions in (2.4): associated cones do not work for $p$
small. This goal will be achieved  in (3.6) where we introduce a
group analogue of ${\mathbb K}({\cal S}(\tilde{h},\tilde{e}))$.
The latter will be linked with $\C(\g)$ by means of a
$G$-equivariant isomorphism between $\cal U$ and $\N$.

In fact, we need the $G$-equivariant isomorphism between ${\cal
U}$ and $\N$ introduced in [\cite{BR}]. This isomorphism, call it
$\eta$,  is defined in loc. cit. as follows: If $G=\GL(V)$, one
just puts $\eta(u)=u-1$ for all $u\in\cal U$. If $G$ is not of
type $\mathrm A$,  then $\g$ is a simple Lie algebra. In this case
one picks a finite dimensinal rational representation
$\rho\colon\,G\rightarrow\GL(W)$ with $\text{ker}\,\rho\subseteq
Z(G)$ such that the trace form $(X,Y)\longmapsto
\,\text{tr}\,(\text{d}\rho)_e(X)\circ(\text{d}\rho)_e(Y)$ on $\g$
is nondegenerate; see [\cite{SpSt}, Chap.~I, Lemma~5.3].  One then
identifies $\g$ with its image $(\text{d}\rho)(\g)$. By the choice
of $\rho$, there is a subspace $M\subset {\frak gl}(W)$ such that
${\frak gl}(W)=\,\g\oplus M$ and $[\g,M]\subseteq M$.
Specifically, $M$ is the orthogonal complement of $\g$ with
respect to the trace form on ${\frak gl}(W)$. Finally, one takes
for  $\eta$ the restriction of $\text{pr}_1\circ \rho$ to $\cal U$
where $\text{pr}_1\colon\, {\frak gl}(W)\twoheadrightarrow \g$ is
the first projection. According to [\cite{BR}, Cor.~9.3.4],
$\eta\colon\,{\cal U}\,\stackrel{\sim}{\longrightarrow}\,\N$ is a
$G$-equivariant isomorphism of algebraic varieties. Hence so is
the map $\eta\times \text{id}_\N\colon\, {\cal
U}\times\N\,\longrightarrow\,\N\times \N$.

We now define ${\frak X}:= \{(u,x)\in {\cal U}\times \N\,|\,
(\mbox{Ad}\, u)\,x=x\},$ a  closed subset in $\cal U\times \N$,
and denote by $\tilde{\eta}$ the restriction of $\eta\times
\mathrm{id}_{\N}$ to $\frak X$.
\begin{lemma} The map $\tilde{\eta}$
is a $G$-equivariant isomorphism between
$\frak X$ and ${\frak C}^{\mathrm nil}(\g)$. For any $e\in\N$
it maps $(R_u(Z_G(e)),e)$ onto
$(\Lie\,R_u(Z_G(e)),e)$.
\end{lemma}
\pf Let $(u,e)\in\frak X$. Then $(\text{Ad}\,u)\,e=e$ forcing
$\rho(u)e\rho(u)^{-1}=e$ (recall that we identify $\g$ with
$(\text{d}\rho)(\g)$). Thus $\rho(u)$ commutes with $e$. Write
$\rho(u)=x+m$ with $x\in\g$ and $m\in M$. Then $x=\eta(u)\in\N$
and $0=[e,\rho(u)]=[e,x]+[e,m].$ Hence $[e,\eta(u)]=0$. This shows
that $\tilde{\eta}(\frak X)\subseteq \C(\g)$.

As a consequence, $\tilde{\eta}$ sends $({\cal U}\cap Z_G(e),e)$
to $(\N\cap\z(e),e)$. From [\cite{Jan}, Sect.~5] and [\cite{P02},
Sect.~2] we know that  ${\cal U}\cap Z_G(e)$ is isomorphic to
${\cal U}\big(Z(\lambda_e)\cap Z_G(e)\big)\times\, R_u(Z_G(e))$.
Our discussion in (3.1) shows that
$$\N\cap\z(e)\,\cong\,\N\big(\Lie\,Z(\lambda_e)\cap
Z_G(e)\big)\times\,\Lie\,R_u(Z_G(e))$$ is an irreducible variety.
Since $Z(\lambda_e)\cap Z_G(e)$ is a reductive group, the
varieties $\N\big(\Lie\,Z(\lambda_e)\cap Z_G(e)\big)$ and
 ${\cal U}\big(Z(\lambda_e)\cap Z_G(e)\big)$ have the same
dimension. This implies that $\dim\,{\cal U}\cap Z_G(e)\,=\,\dim\,\N\cap\z(e)$.

Since $\eta\times\text{id}_\N$ is an isomorphism of varieties, it maps the closed
set $({\cal U}\cap Z_G(e),e)$ onto a closed subset of
$(\N\cap\z(e),e)$ of the same dimension. The irreducibility of $(\N\cap\z(e),e)$ now
yields that $\tilde{\eta}$ maps $({\cal U}\cap Z_G(e),e)$ onto
$(\N\cap \z(e),e)$. Since this holds for any $e\in\N$, it must be that
$\tilde{\eta}({\frak X})\,=\,\C(\g)$. But
then $\tilde{\eta}\colon\,{\frak X}\,\longrightarrow\,
\C(\g)$ is a $G$-equivariant isomorphism of varieties.

Let $P(\lambda_e)$ be the parabolic subgroup of $G$ with Lie
algebra ${\frak p}(\lambda_e):=\,\bigoplus_{i\ge
0}\,\g(i;\lambda_e)$. It follows from [\cite{Mc}, Lemma~28] that
$\eta$ maps $R_u(P(\lambda_e))$ onto $\Lie\,R_u(P(\lambda_e))$.
Note that
$\Lie\,R_u(P(\lambda_e))=\,\bigoplus_{i>0}\,\g(i;\lambda_e)$.
Since $\tilde{\eta}$ maps $({\cal U}\cap Z_G(e),e)$ onto
$(\N\cap\z(e),e)$,  it maps $\big(R_u(P(\lambda_e))\cap
Z_G(e),e\big)$ onto $\big(\Lie\,R_u(P(\lambda_e))\cap
\z(e),e\big)$. Combined with [\cite{Jan}, Sect.~5] or [\cite{P02},
Sect.~2] this yields that $\tilde{\eta}$ maps $(R_u(Z_G(e)),e)$
onto $(\Lie\,R_u(Z_G(e)),e)$ as desired. \qed

\smallskip

\noindent {\bf Remark.} Let ${\frak C}^{\,\rm
unip}(G):=\{(x,y)\in{\cal U}\times{\cal U}\,|\,\,xy\,=\,yx\}$ be
the unipotent commuting variety of $G$. Since both $\rho(\g)$ and
$M$ are stable under all ${\rm Int}\,\rho(g)$ with $g\in G$, the
argument used in the proof of Lemma~3.5 also shows the map
$\text{id}_{\cal U}\times \eta$ induces a $G$-equivariant
isomorphism between ${\frak C}^{\,{\rm unip}}(G)$ and $\frak X$.
As a result, the varieties ${\frak C}^{\,{\rm unip}}(G)$ and
$\C(\g)$ are $G$-equivariantly isomorphic.

\smallskip

\noindent {\bf 3.6.} The group $G$ acts on $G\times\N$ via
$g\cdot(x,n)=((\text{Int}\,g)\,x,({\text Ad}\,g)\,n)$ fot all
$g,x\in G$ and $n\in\N$. Given a distinguished nilpotent element
$e\in\g$ we denote by ${\cal Y}(\lambda_e,e)$ the Zariski closure
of $G\cdot (\lambda_e(k^*),e)$ in $G\times \N$. By the definition
of ${\cal Y}(\lambda_e,e)$, the morphism $$\nu\colon\,G\times
k^*\,\longrightarrow\,{\cal Y}(\lambda_e,e),\qquad\, \,
(g,t)\,\longmapsto \,
\big(g\,\lambda_e(t)\,g^{-1},(\text{Ad}\,g)\,e\big),$$ is
dominant. It follows that ${\cal Y}(\lambda_e,e)$ is an
irreducible variety. Since $k$ is infinite, there is $t_0\in k^*$
such that $Z_G(\lambda_e)\,=\,Z_G\big(\lambda_e(t_0)\big)$. The
fibre $\nu^{-1}(\nu(1,t_0))=\nu^{-1}(\lambda_e(t_0),e)$ is nothing
but $\,\text{Stab}_G(\lambda_e(t_0),e)$, the stabiliser of
$(\lambda_e(t_0),e)$ in $G$. In view of [\cite{Bo}, Chap.~III,
(9.1), (9,4)] and the choice of $t_0$, we have
$$\Lie\big(\text{Stab}_G(\lambda_e(t_0),e)\big)\subseteq \,
\z(e)\cap\Lie\,Z_G\big(\lambda_e(t_0)\big)\,=\,\z(e)\cap
\g(0;\lambda_e)\,=\,\z(e;0).$$

Since $e$ is distinguished, the reductive part of
$Z_G(e)^\circ\,=\,(Z(G)\cdot Z_{(G,G)}(e))^\circ$ coincides with
$Z(G)^\circ$, the connected centre of $G$. Combining this with
[\cite{Jan}, Sect.~5] and [\cite{P02}\, Sect.~2] one derives that
$\z(e;0)$ coincides with $\z(\g)$, the centre of $\g$. By our
assumption on $G$, we have $\z(\g)=\{0\}$ unless $G$ is of type
$\mathrm A,$ in which case $\dim\,\z(\g)=\dim\,Z(G)^\circ=1$. This
implies that $\big(\text{Stab}_G(\lambda_e(t_0),e)\big)^\circ
=\,Z(G)^\circ$. Consequently,
$\dim\,\nu^{-1}(\nu(1,t_0))=\dim\,Z(G)^\circ$. On the other hand,
it is easy to see that $\dim\,\nu^{-1}(\nu(g,t))\ge
\dim\,Z(G)^\circ$ for all $(g,t)\in G\times k^*$. Applying the
theorem on the dimensions of the fibres of a morphism we now get
\begin{eqnarray}
\dim\,{\cal Y}(\lambda_e,e)=\dim\,G-\dim\,Z(G)^\circ+1=\dim\,(G,G)+1.
\end{eqnarray}

Now let $l=\text{rk}\,(G,G)$ and let $\rho_1,\ldots,\rho_l$ be the
fundamental representations of $G$ (when $G= GL(V)$, we take as
$\rho_i$ the $i$th exterior power of the vector representation of
$G$). For $1\le i\le l$, define the regular invariant function
$\chi_i$ on $G$ by setting $\chi_i(g)=\text{tr}\,\rho_i(g)$ for
all $g\in G$. Define $\chi_0\in k[G]^G$ by setting, for all $g\in
G$, $$\chi_0(g)\,=\,\left\{
\begin{array}{ll}
\det\,g\,&\mbox{  when $G$ is of type }{\mathrm A},\\ 1\, &\mbox{
otherwise}.\end{array}\right.$$ Define the morphism $F\colon{\cal
Y}(\lambda_e,e)\,\longrightarrow\, {\mathbb A}^{l+1}$ by letting
$$F(g,n)=\big(\chi_0(g),\chi_1(g),\ldots,\chi_l(g)\big)\qquad
\,\forall\,(g,n)\in{\cal Y}(\lambda_e,e),$$ and denote by ${\cal
Y}_0(\lambda_e,e)$ the fibre $F^{-1}(F(1,0))$.
\begin{prop}
${\cal Y}_0(\lambda_e,e)$ is a closed subset in ${\frak X}$ and
all irreducible components of ${\cal Y}_0(\lambda_e,e)$ have
dimension equal to $\dim\,(G,G)$.
\end{prop}
\begin{pf}
(1) First note that $(G,G)=\{g\in G\,|\,\chi_0(g)=1\}$. Together
with [\cite{St}, (3.4), Cor.~4] this implies that ${\cal U}=\{g\in
G\,|\,\chi_i(g)=\chi_i(1)\,\,\mbox{ for all }i\}$ (one should take
into account that the group $(G,G)$ is simply connected). It
follows that ${\cal Y}_0(\lambda_e,e)\subseteq {\cal U}\times \N.$
On the other hand, the set ${\cal R}:=\{(g,x)\in G\times\N\,|\,
(\text{Ad}\,g)\,x\in kx\}$ is Zariski closed and $G$-stable in
$G\times \N$. As $(\lambda_e(k^*),e)\subset \cal R$, this yields
${\cal Y}(\lambda_e,e)\subseteq {\cal R}$. If $u\in\cal U$, then
$1$ is the only eigenvalue of $\text{Ad}\,u$. Hence $({\cal
U}\times\N)\cap{\cal R}\subseteq {\frak X}$ forcing $${\cal
Y}_0(\lambda_e,e)\subseteq ({\cal U}\times\N)\cap{\cal
R}\subseteq\frak X.$$

\smallskip

\noindent (2) Let $\cal C$ denote the Zariski closure of
$\{F(\lambda_e(t),0)\,|\,t\in k^*\}$ in ${\mathbb A}^{l+1}$. Since
$\lambda_e\colon\, k^*\rightarrow G$ is a nontrivial homomorphism,
$\cal C$ is an affine curve. Since all $\chi_i$ are $G$-invariant,
$F$ maps ${\cal Y}(\lambda_e,e)$ onto a Zariski open subset of
$\cal C$, say ${\cal C}_0$. Clearly, ${\cal C}_0$ is a curve on
its own and the morphism $F\colon\,{\cal
Y}(\lambda_e,e)\,\longrightarrow\,{\cal C}_0$ is surjective. As
${\cal Y}(\lambda_e,e)$ is irreducible, we now combine [\cite{Sh},
Chap.~I, \S~6, Theorem~7] with (11) to deduce that all irreducible
components of $F^{-1}(F(1,0))= {\cal Y}_0(\lambda_e,e)$ have
dimension equal to $\dim\,(G,G)$. This completes the proof.
\end{pf}

\noindent {\bf 3.7.} Now we are ready for the main result of this
section.
\begin{thm}
Let $G$ be a connected reductive group over $k$ and assume that
$p={\mathrm char}\,k$ is a good prime for $G$. Let
$\g=\Lie\,G$ and $n=\dim\,(G,G)$.
Then all irreducible components of
${\frak C}^{\mathrm nil}(\g)$ are $n$-dimensional and their number
equals
the number of distinguished nilpotent conjugacy classes in $\g$.
Moreover, each irreducible component of
${\frak C}^{\mathrm nil}(\g)$ has the form ${\frak C}(e)$ for some
distinguished nilpotent element $e\in\g$.
\end{thm}
\pf (1) As explained in (3.1), it can be assumed in proving this
theorem that $G$ is either $\GL(n)$ or a simple algebraic group of
type different from $\mathrm A$. Let $e$ be an almost
distinguished nilpotent element in $\g$. Thanks to Proposition~3.3
it suffices to show that ${\frak C}(e)\subseteq {\frak C}(e')$ for
some distinguished $e'\in\N$. Since ${\frak C}(e)\,=\,{\frak
C}((\text{Ad}\,g)\,e)$ for any $g\in G$ we may assume that $e$ is
as in (3.4); in particular, $e$ is Richardson in ${\frak u}_{I,J}$
for some $(I,J)\in{\cal P}(\Pi)$. We now pick a distinguished
$\tilde{e}\in\N$ following the procedure described in (3.4) and
take as $\lambda_{\tilde{e}}$ the $1$-parameter subgroup
$\lambda_J$.

\smallskip

\noindent (2) We claim that
$\big(\lambda_{\tilde{e}}(k^*)\,R_u(Z_G(e)),e\big)\subset\, {\cal
Y}(\lambda_{\tilde{e}},\tilde{e})$. To that end we let
$b_1,\ldots, b_m$ be the weights of
$(\text{Ad}\,\lambda_{\tilde{e}})(k^*)$ on $\Lie R_u(Z_G(e))$. By
Lemma~3.4, none of them is zero. It follows that
$k^*_{\scriptstyle{\text{reg}}}:=\{t\in
k^*\,|\,b_i(\text{Ad}\,\lambda_{\tilde{e}}(t))\ne 1 \ \, \mbox{for
all } \,i\}$ is a nonempty Zariski open subset in $k^*$. We pick
any $v\in k^*_{\scriptstyle{\text{reg}}}$ and set
$s=\lambda_{\tilde{e}}(v)$, $\tilde{\cal Y}={\cal
Y}(\lambda_{\tilde{e}},\tilde{e})$,
$L\,=\,Z_{G}(\lambda_{\tilde{e}})$. By [\cite{P02}, Sect.~2], the
orbit $(\text{Ad}\,L)\,\tilde{e}$ is open in
$\g(2;\lambda_{\tilde{e}})$. Since $\tilde{\cal Y}$ is $L$-stable
and Zariski closed, we then have
$(s,\g(2;\lambda_{\tilde{e}}))\subset\, \tilde{\cal Y}.$ Since
$e\in\g(2;\lambda_{\tilde{e}})$ by our discussion in (3.4), we
deduce that $(s,e)\in\tilde{\cal Y}$.

Let $U(s)$ denote the centraliser of $s$ in $R_u(Z_G(e))$. As $s$
is semisimple and $R_u(Z_G(e))$ is a connected unipotent group,
$U(s)$ is connected too; see [\cite{Bo}, Chap.~III, (9.3)]. Since
$\Lie\,U(s)=\big(\Lie\,R_u(Z_G(e))\big)^s$, by
[\cite{Bo},Chap.~III, (9.4)], and since $b_i(\text{Ad}\,s)\ne 1$
for all $i$, it must be that $U(s)=\{1\}$. As $s$ normalises the
unipotent group $R_u(Z_G(e))$, the orbit $
\big(\text{Int}\,R_u(Z_G(e))\big)s$ is Zariski closed in
$s\,R_u(Z_G(e))$. By the preceding remark, it has the same
dimension as $s\,R_u(Z_G(e))$. The irreducibility of
$s\,R_u(Z_G(e))$ now yields $(s\,R_u(Z_G(e)),e)\subset\tilde{\cal
Y}$. Since this holds for all $v\in
k^*_{\scriptstyle{\text{reg}}}$, we get
$(\lambda_{\tilde{e}}(k^*_{\scriptstyle{\text{reg}}})\,y,e)\subset\tilde{\cal
Y}$ for all $y\in R_u(Z_G(e))$. But then
$(\lambda_{\tilde{e}}(k^*)\,y,e)=\,
\overline{(\lambda_{\tilde{e}}(k^*_{\scriptstyle{\text{reg}}})\,y,e)}
\subset\,\tilde{\cal Y}$ for all $y\in R_u(Z_G(e)),$ and the claim
follows.

\smallskip

\noindent
(3) Proposition~3.6 now shows that $\big(R_u(Z_G(e)),e\big)\subset\,
\tilde{\cal Y}\cap {\frak X}\,=\,
{\cal Y}_0(\lambda_{\tilde{e}},\tilde{e}).$ Combined with Lemma~3.5 this
implies that
$(\Lie\,R_u(Z_G(e)),e)\subset \tilde{\eta}\big({\cal Y}_0(\lambda_{\tilde{e}},\tilde{e})
\big)\subseteq \C(\g)$  and all irreducible components of the Zariski closed set
$\tilde{\eta}\big({\cal Y}_0(\lambda_{\tilde{e}},\tilde{e})\big)$
have dimension $n$.

Let $Z$ be an irreducible component of $\tilde{\eta}
\big({\cal Y}_0(\lambda_{\tilde{e}},\tilde{e})\big)$ containing
$(\Lie\,R_u(Z_G(e)),e)$.
Then $Z\subseteq {\frak C}(e')$ for some $e'\in\N$; see Proposition~3.3.
By our discussion in (3.1), the variety $\N(\z(e'))=\,\N\cap\z(e')$ is irreducible.
So it is immediate from the
definition of ${\frak C}(e')$ that the morphism $$\xi\colon\,G\times
\N(\z(e'))\,\longrightarrow\, {\frak C}(e'),\quad\ \ \,
\xi(g,x)\,=\, \big((\text{Ad}\,g)\,e',(\text{Ad}\,g)\,x\big),$$
is dominant. For every $x\in\N(\z(e'))$ the fibre
$\xi^{-1}(\xi(1,x))$ is nothing but the set of all pairs
$\big(g,(\text{Ad}\,g)^{-1}\,x\big)$ with $g\in Z_{G}(e')$. It
follows that $\xi^{-1}(\xi(1))\,\cong\, Z_{G}(e')$ as varieties.
The theorem on the dimension of the fibres of a morphism
in conjunction with our discussion in (3.1) now gives
\begin{eqnarray*}
\dim\,{\frak
C}(e')&=&\dim\,G+\dim\,\N(\z(e'))-\dim\,Z_{G}(e')\\
&=&\dim\,G+\dim\N(\z(e'))-\dim\,\z(e')\\
&=&\dim\,G-\text{rk}\big(Z(\lambda_{e'})\cap Z_G(e')\big)\\
&=&\dim\,(G,G)-\text{rk}\big(Z(\lambda_{e'})\cap Z_{(G,G)}(e')\big).
\end{eqnarray*}
As $\dim\,Z=n$, we get $\text{rk}\big(Z(\lambda_{e'})\cap
Z_{(G,G)}(e')\big)=0$. As $Z_G(e')\,=\,Z(G)\cdot Z_{(G,G)}(e')$,
our discussion in (3.1) shows that the group $Z_{(G,G)}(e')^\circ$
must be unipotent.

As a result, $e'$ is distinguished in $\g$ and $Z=\,{\frak
C}(e')$, by dimension reasons. Since $e$ is almost distinguished,
our remark in (3.1) yields that $\N(\z(e;0))=\{0\}$ and
$\N\cap\z(e)=\,\textstyle{\bigoplus_{i>0}}\,\z(e;i)\,=\,\Lie\,R_u(Z_G(e)).$
As ${\frak C}(e')$ is $G$-stable, we then have $G\cdot (\N\cap
\z(e),e)\,\subseteq \,{\frak C}(e')$ forcing ${\frak
C}(e)\subseteq {\frak C}(e')$. This completes the proof. \qed

\section{\bf Applications and concluding remarks}

\noindent {\bf 4.1.} Let $S$ be a smooth projective surface over
$k$ and let $\text{Hilb}^n(S)$ be the Hilbert scheme parametrising
the $0$-dimensional subschemes of $S$ of length $n$. The scheme
$\text{Hilb}^n(S)$ is known to be smooth and projective; see
[\cite{F1}, \cite{F2}]. Let $\text{Sym}^n(S)$ denote the $n$-fold
symmetric power of $S$, the geometric quotient of $S^n$ by the
permutation action of ${\frak S}_n$. By assigning to any
$0$-dimensional subscheme of $S$ its support with multiplicities
one obtains a map
$$w_n\colon\,\text{Hilb}^n(S)\,\longrightarrow\,\text{Sym}^n(S)$$
called the {\it Hilbert-Chow morphism}. According to [\cite{F1},
\cite{F2}], this map provides a natural desingularisation of
$\text{Sym}^n(S)$. If $z$ is a point of $\text{Sym}^n(S)$
representing the effective $0$-cycle  $\sum n_i\,[s_i]$ with $\sum
n_i=n$, then the fibre $w_n^{-1}(z)$ is isomorphic to
$\prod_i\text{Hilb}^{n_i}\,{\cal O}_{S,\,s_i}$ where
$\text{Hilb}^r\,{\cal O}_{S,s}$ is the punctual Hilbert scheme
parametrising the ideals of colength $r$ in the local ring ${\cal
O}_{S,s}$.

For any $s\in S$ the scheme $(\text{Hilb}^r\,{\cal
O}_{S,s})_{\text{red}}$ is isomorphic to the (reduced) Hilbert
scheme ${\cal H}_r=\text{Hilb}^r\,k[[x,y]]$ parametrising the
ideals of colength $d$ in the ring of formal power series
$k[\!\,[x,y]]$; see [\cite{I}, Chap.~1]. This shows that the
punctual Hilbert schemes ${\cal H}_d$ with $r\le n$ serve as the
building blocks for {\it all} fibres of the Hilbert-Chow morphism.

Let $\frak m$ be the maximal ideal of $k[[x,y]]$ and
$A=k[[x,y]]/{\frak m}^r$, a finite dimensional local algebra over
$k$. Let $\text{Gr}^r(A)$ denote the Grassmannian of all
$k$-subspaces of codimension $r$ in $A$. By [\cite{I}, Chap.~1],
any ideal $I$ of colength $r$ in $k[[x,y]]$ contains ${\frak m}^r$
and the map $I\mapsto I/{\frak m}^r$ identifies ${\cal H}_r$ with
a closed subset of $\text{Gr}^r(A)$. The algebraic variety
structure on ${\cal H}_r$ is given by this identification.

Now let $V$ is an $r$-dimensional vector space over $k$ and
$\g={\frak gl}(V)$. Observe that $G=\GL(V)$ acts on $\C(\g)\times
V$ by $g\cdot (a,b;v)=(gag^{-1},gbg^{-1};g(v))$ for all $g\in G$
and all $(a,b;v)\in \C(\g)\times V$. Let $U$ be the Zariski open
subset in $\C(\g)\times V$ consisting of all $(a,b;v)$ such that
$a^ib^j(v)$ with $i,j\in\mathbb Z_+$ span $V$. Clearly, $U$ is
$G$-stable. As observed in [\cite{Na}, \cite{Bar}], the morphism
$\psi\colon\,U\longrightarrow\, \text{Gr}^r(A)$ assigning to
$(a,b;v)\in U$ the point $\{\phi\in A\,|\,\phi(a,b)(v)=0\}$ on
$\text{Gr}^r(A)$, maps $U$ onto ${\cal H}_r$. Moreover, $\GL(V)$
acts freely on $U$ and the fibres of $\psi$ are exactly the
$\GL(V)$-orbits on $U$.

It is immediate from Theorems~2.5 and 3.7 that $\C(\g)\times V$ is
an irreducible variety of dimension $r^2+r-1$. Then so is $U$
being open dense in $\C(\g)\times V$. Combining the theorem on the
dimension of the fibres of a morphism with the discussion above we
now obtain the following:
\begin{cor}
If $k$ is an algebraically closed field of characteristic $p\ge
0$, then the punctual Hilbert scheme ${\cal H}_r$ is irreducible
over $k$ and $\dim \,{\cal H}_r=\,r-1$.
\end{cor}
This corollary generalises earlier results by Brian\c{c}on
[\cite{Br}] for $p=0$, Iarrobino [\cite{I}] for $p>n$,  and Basili
[\cite{Bas}] for $p\ge n/2$. In conjunction with an observation in
[\cite{Na}, (1.2)] it also shows that ${\cal H}_r$ is isomorphic
(as a reduced scheme) to the closure in $\text{Gr}^r(A)$ of the
$(r-1)$-dimensional family of  principal ideals $${\cal
J}=\big\{A\big(x-(t_1y+t_2y^2+\cdots+t_{r-1}y^{r-1})\big)\,|\
(t_1,t_2,\ldots,t_{r-1})\in {\mathbb A}^{r-1}\big\}.$$

\smallskip

\noindent
{\bf 4.2.} In this subsection we assume that $k$ is an algebraically closed field
of characteristic $p>0$. Let $\L$ be a finite dimensional restricted Lie algebra
over $k$ with $p$th power map $x\mapsto\, x^{[p]}$. We denote by $\N(\L)$ the
nilpotent variety of $\L$, the set of all $x\in\L$ with $x^{[p]^N}=0$ for $N\gg 0$.
It is immediate from Jacobson's formula [\cite{Jac}, Chap.~V, Sect.~7] that
the map
$\pi\colon\,\L\rightarrow\L,\ \,x\mapsto x^{[p]},$ is a morphism
given by homogeneous polynomial functions
on $\L$ of degree $p$.

Recall that an element $x\in\L$ is {\it semisimple} if it lies in
the restricted subalgebra of $\L$ generated by $x^{[p]}$. Let
$\L_{ss}$ denote the set of all semisimple elements of $\L$. Using
the Jordan-Chevalley decomposition in $\L$ it is easy to observe
that $\L_{ss}$ coincides with $\pi^N(\L)$ for $N\gg 0$. We define
$$e(\L):=\min\,\{r\in{\mathbb Z}_+\,|\,\pi^r(W)\subseteq \L_{ss}
\text{ for some open dense }W\subseteq \L\}.$$

Recall that a restricted subalgebra $\frak t\subseteq \L$ is
called {\it toral} (or a {\it torus}) if the $p$th power map of
$\L$ is injective on $\frak t$. By [\cite{Jac}, Chap.~V, Sect. 8],
any torus $\frak t$ in $\L$ is  abelian and spanned by its subset
${\frak t}^{\scriptstyle{\text{tor}}}:=\,\{t\in{\frak t}\,
|\,t^{[p]}=t\}$. We define $$MT(\L):=\max\,\{\dim\,{\frak
t}\,|\,\,{\frak t} \text{ is toral in }\L\}.$$

Note that $MT(\L)=0$ if and only if $\L=\N(\L)$, while
from the main results of [\cite{P86}] it follows that $e(\L)=0$ if and
only if $\L$ possesses a toral Cartan subalgebra (see also
[\cite{P89}, Theorem~1]). Let $n=\dim\,\L,$ $e=e(\L),$ $s=MT(\L)$, and let
$v_1,\ldots,v_n$ be a basis of $\L$. According to [\cite{P89}, Theorem~2],
there exist homogeneous polynomials $\psi_0,\ldots,\psi_{s-1}\in
k[X_1,\ldots,X_{n}]$ with $\deg\,\psi_i=p^{s+e}-p^{i+e}$ and the property
that for any
$x=\sum_{i=1}^n x_i\,v_i$   one has
\begin{eqnarray}x^{[p]^{s+e}}=\,\,
\textstyle
{\sum_{i=0}^{s-1}}\,\psi_i(x)\,
x^{[p]^{i+e}}.\end{eqnarray}
We view
$\psi_0,\ldots,\psi_{s-1}$ as polynomial functions  on $\L$.
\begin{thm}
The following are true:
\begin{itemize}

\item[1.] $\psi_i(x^{[p]})=\psi_i(x)^p$ for all $x\in\L$ and all $i\le s-1$.

\smallskip
\item[2.] $\N(\L)$ coincides with the set of all common zeros of
$\psi_0,\ldots, \psi_{s-1}$ in $\L$.

\smallskip
\item[3.] All irreducible components of $\N(\L)$ have dimension $n-s$.
\end{itemize}
\end{thm}
\pf (1) Let $g_{i+1}=\psi_{i}^p-\psi_{i}\circ\pi$ where $0\le i\le
s-1$. Then $g_{i+1}$ is a polynomial function on $\L$. Applying
$\pi$ to both sides of (12) we get $x^{[p]^{s+e+1}}\!=\,\,
\textstyle{\sum_{i=0}^{s-1}}\,\psi_i(x)^p\,x^{[p]^{i+e+1}}$ while
substituting in (12) $x$ by $x^{[p]}$ we get
$x^{[p]^{s+e+1}}\!=\,\,
\textstyle{\sum_{i=0}^{s-1}}\,\psi_i(x^{[p]})\,x^{[p]^{i+e+1}}.$
Hence
\begin{eqnarray}
\textstyle{\sum_{i=1}^{s}}\,g_i(x)\,x^{[p]^{e+i}}=\,0\qquad \quad\,  (\forall\,x\in\L).
\end{eqnarray}
For $r\in\mathbb N$ and $x\in\L$, write
$\pi^r(x)=\sum_{i=1}^n\,\pi_{r,i}(x)\,v_i$ where $\pi_{r,i}$ is a
homogeneous polynomial function on $\L$ of degree $p^r$. As
explained in [\cite{P89}, \S~1] the matrix
$$\left(\begin{array}{cccc} \pi_{e,1}&\pi_{e,2}&\cdots
&\pi_{e,n}\\ \pi_{e+1,1}&\pi_{e+1,2}&\cdots &\pi_{e+1,n}\\
\vdots&\vdots & &\vdots\\ \pi_{e+s-1,1}&\pi_{e+s-1,2}&\cdots &
\pi_{e+s-1,n}
\end{array}\right)$$ has rank $s$.
It follows that there is an open dense subset $U\subset \L$ with
the property that $u^{[p]^e}, u^{[p]^{e+1}}, \ldots, \,
u^{[p]^{e+s-1}}$ are linearly independent for all $u\in U$. Let
$W$ be an open dense set in $\L$ with $\pi^e(W)\subset\L_{ss}$.
Note that for each $u\in W$ the linear span of $u^{[p]^j}$ with
$j\ge e$ is a torus. It follows that $u^{[p]^{e+1}},
u^{[p]^{e+2}}, \ldots, \, u^{[p]^{e+s}}$ are linearly independent
whenever $u\in U\cap W$. So (13) shows that all $g_i$ vanish on
$U\cap W$. As $U\cap W$ is open dense in $\L$, we get
$\psi_i\circ\pi=\psi_i^p$ for all $i\le s-1$.

\smallskip

\noindent (2) Let ${\cal Z}$ denote the set of all common zeros of
$\psi_0,\ldots,\psi_{s-1}$ in $\L$. It is immediate from (12) that
$\N(\L)\supseteq \cal Z$. On the other hand, if $x\in \N(\L)$,
then $$\psi_i(x)^{p^N}=\psi_i(x^{[p]^N})=\psi_i(0)=0 \qquad\quad\
(0\le i\le s-1)$$ provided $N\gg 0$. Hence $\cal Z=\N(\L)$
implying that each irreducible component of $\N(\L)$ has dimension
$\ge n-s$, by the affine dimension theorem. If $\N(\L)$ had a
component of dimension $>n-s$, then it would intersect
nontrivially with any $s$-dimensional torus in  $\L$ (again by the
affine dimension theorem). But this is impossible as no nonzero
element in $\L$ can be both nilpotent and semisimple. Thus all
irreducible components of $\N(\L)$ are of dimension $n-s$. \qed

\smallskip

\noindent {\bf Remark.} It is conjectured in [\cite{P89}] that the
variety $\N(\L)$ is always irreducible. This conjecture is still
open; it is not clear at present how to approach it without using
the structure theory of finite dimensional restricted Lie
algebras.

\smallskip

\noindent {\bf 4.3.} We are going to apply Theorem~4.2 for
estimating the dimension of $\C(\g)$ in a more general setting. We
assume in this subsection that $k$ is an algebraically closed
field of characteristic $p>0$ and $G$ is a connected reductive
$k$-group (so we allow $p$ to be a bad prime from now on). We
adopt the notation introduced in Sections~2 and 3.

Let $e_1,\ldots, e_m$ be representatives of the nilpotent
$G$-orbits in $\g$; these orbits are finite in number, by
[\cite{HSp}]. Then $\C(\g)= {\frak C}(e_1)\cup\ldots\cup{\frak
C}(e_m)$, as before, but in the present case there is no guaranty
that the ${\frak C}(e_i)$'s are all irreducible. However, the
nilpotent variety $\N$ of $\g$ is still irreducible, by
[\cite{Bo}, Chap.~IV, (14.16)]. Combined with [\cite{HSp}] this
implies that $\N$ contains a unique open $G$-orbit. Let $e_{\rm
reg}$ be a representative of this orbit.
\begin{prop} The following are true:
\begin{itemize}
\item[1.] For any $e\in \N$ the variety ${\frak C}(e)$ is equidimensional and
$$\dim {\frak C}(e)\,=\,\dim G-\dim Z_G(e)+\dim \z(e)-MT(\z(e)).$$

\smallskip

\item[2.] ${\frak C}(e_{\rm reg})$ is an irreducible component of
${\frak C}^{\rm nil}(\g)$ and $$\dim {\frak C}(e_{\rm reg})\,=\,\dim G-{\rm rk}\,G+
\dim \z(e_{\rm reg})-\dim \z(\g).$$
\end{itemize}
\end{prop}
\pf (1) Let $M_1,\ldots, M_s$ be the irreducible components of
$\N(\z(e))$. By Theorem~4.2(3), $\dim M_i=\dim \z(e)-MT(\z(e))$
for all $i\le s$. Let $\tilde{M}_i$ denote the Zariski closure of
$G\cdot (e,M_i)$ in $\g\times \g$. Then ${\frak
C}(e)=\tilde{M}_1\cup\ldots\cup \tilde{M}_s$. It follows that each
irreducible component of ${\frak C}(e)$ has the form $\tilde{M}_j$
for some $j\le s$. The map $$\xi\colon\,G\times
M_j\,\longrightarrow \,\tilde{M}_j, \quad \ \ \,
(x,m)\longmapsto\,\big((\text{Ad}\,x)\,e,(\text{Ad}\,x)\,m)\big),$$
is a dominant morphism of algebraic varieties. The fibre
$\xi^{-1}(\xi(g,a))$ consists of all $(x,m)\in G\times M_j$ with
$g^{-1}x\in Z_G(e)$ and $(\text{Ad}\,x^{-1}g)\,a=m$. From this it
is immediate that $\xi^{-1}(\xi(g,a))\cong \text{Tran}\,
(a,M_j):=\{z\in Z_G(e)\,|\, \,(\text{Ad}\, z)\,a\in M_j\}.$ Since
$Z_G(e)$ permutes the components $M_1,\ldots, M_s$, the connected
component $Z_G(e)^\circ$ acts on $\text{Tran}\,(a,M_j)$ by left
translations. Since $Z_G(e)^\circ$ has finite index in $Z_G(e)$,
it follows that $\text{Tran}\,(a,M_j)$ and $Z_G(e)$ have the same
dimension. As a result, $\dim \xi^{-1}(\xi(g,a))=\dim Z_G(e)$. The
theorem on the dimension of the fibres of a morphism now gives
\begin{eqnarray*}
\dim \tilde{M}_j&=&\dim G+\dim M_j-\dim Z_G(e)\\
&=&\dim G-\dim Z_G(e)+\dim\z(e)-MT(\z(e)),\end{eqnarray*}
proving the first part of the proposition. As
$\dim Z_G(e_{ \rm reg})=
\text{rk}\, G$, we also obtain \begin{eqnarray}\dim {\frak C}(e_{\rm reg})=\,
\dim G-\text{rk}\,G+\dim \z(e_{\rm reg})-MT(\z(e_{\rm reg})).\end{eqnarray}

\smallskip

\noindent (2) Let $T$ be a maximal torus in $G$ and
$\t\,=\,\Lie\,T$. Let $\{\alpha_1,\ldots,\alpha_l\}$ be a basis of
simple roots in the root system of $G$ relative to $T$. According
to [\cite{Jan}, Prop.~4.14], it can be assumed that $e_{\rm
reg}=e_{\alpha_1}+\ldots + e_{\alpha_l}$ where $e_{\alpha_i}$ are
simple root vectors in $\g$ (the argument in [\cite{Jan}] is based
on [\cite{Ca}, Prop.~5.8.5]) which is applicable in any
characteristic). There is a cocharacter $\lambda_{\rm reg}\in
X_*(T)$ such that $\alpha_{i}(\lambda_{\rm reg}(t))=t^2$ for all
$i$ and all $t\in k^*$. It is easy to see that $\g(0;\lambda_{\rm
reg})=\,\t$. For $r\in\mathbb Z,$ we  set $\z(e_{\rm
reg};\,r)\,=\,\z(e_{\rm reg})\cap\g(r;\lambda_{\rm reg})$. Then
$\z(e_{\rm reg};\,0)$ consists of all $x\in \frak t$ with
$(\text{d}\alpha_i)_{e}(x)=0$ for all $i$. Since any root $\gamma$
is a $\mathbb Z$-linear combination of simple roots, the root
function $(\text{d}\gamma)_e \in{\frak t}^*$ is an ${\mathbb
F}_p$-linear combination of the root functions
$(\text{d}\alpha_i)_e$. It follows that any root function
$(\text{d}\gamma)_e$ vanishes on $\z(e_{\rm reg};\,0)\subseteq
\frak t$. Since $\frak t$ is abelian, this implies that $\z(e_{\rm
reg};\,0)\subseteq\,\z(\g)$. On the other hand, decomposing
$\z(\g)$ into weight spaces with respect to $T$ one observes
without difficulty that $\z(\g)\subseteq \frak t$. This yields
$\z(e_{\rm reg};\,0)=\,\z(\g)$.

As a result, for any $r\in\mathbb Z$ the set ${\frak
Z}_r:=\{\text{ad}\,x\,|\,x\in\z(e_{\rm reg};\,r)\}$ consists of
nilpotent endomorphisms of $\g$. Since $[{\frak Z}_i,{\frak
Z}_j]\subseteq {\frak Z}_{i+j}$ for all $i,j\in\mathbb Z$, the
union $\bigcup_{\,r\in\mathbb Z}\,{\frak Z}_r$ is a weakly closed
nilset in ${\rm End}\,\g$. Applying Jacobson's theorem on weakly
closed nilsets, see [\cite{Jac}, Chap.~II, \S~2, Theorem~1], we
now obtain that the associative subalgebra generated by the
${\frak Z}_i$'s is nilpotent. But then $\z(e_{\rm reg})$ is a
(restricted) nilpotent subalgebra of $\g$, by Engel's theorem. It
follows that the set $\z(e_{\rm reg})_{ss}$ is central and,
moreover, a unique maximal torus in $\z(e_{\rm reg})$. Since
$\z(e_{\rm reg})_{ss}$ is stable under $\text{Ad}\,\lambda_{\rm
reg}(k^*)$, it must be that $\z(e_{\rm reg})_{ss}=\,\z(\g)$. As a
consequence,
\begin{eqnarray}MT(\z(e_{\rm reg}))=\dim \z(\g).\end{eqnarray}
\noindent (3) Thanks to (14) and (15) we are now reduced to show
that ${\frak C}(e_{\rm reg})$ is a component of $\C(\g)$. Let us
first show that ${\frak C}(e_{\rm reg})$ is irreducible. Let
$\overline{\z}=\z(e_{\rm reg})/\z(\g)$ and let
$\beta\colon\,\z(e_{\rm reg})\twoheadrightarrow \overline{\z}$
denote the canonical homomorphism.Using the Jordan-Chevalley
decomposition in $\z(e_{\rm reg})$ and our final remarks in part~2
one easily observes that $\beta$ induces a bijective morphism
$\N(\z(e_{\rm reg}))\rightarrow \overline{\z}$. The latter, in
turn, induces a bijective morphism ${\mathbb P}\big(\N(\z(e_{\rm
reg}))\big)\longrightarrow\, {\mathbb P}(\overline{\z})$ of
projective varieties. Since ${\mathbb P}(\overline{\z})$ is
obviously irreducible, so must be ${\mathbb P}\big(\N(\z(e_{\rm
reg}))\big)$, by [\cite{Sh}, Chap.~I, \S~6, Theorem~8]. But then
$\N(\z(e_{\rm reg}))$ is irreducible as well, and hence so is
${\frak C}(e_{\rm reg})$.

Let $Z$ be an irreducible component of $\C(\g)$ containing ${\frak
C}(e_{\rm reg})$. By part~1, $Z=\overline{G\cdot(e,M_j)}$ where
$e\in\N$ and $M_j$ is a component on $\N(\z(e))$. It follows that
the first projection $\text{pr}_1\colon\,\g\times
\g\rightarrow\,\g,\ (x,y)\mapsto x$ takes $Z$ to the Zariski
closure of the $G$-orbit of $e$. On the other hand,
$\text{pr}_1({\frak C}(e_{\rm reg}))$ contains the $G$-orbit of
$e_{\rm reg}$, an open subset in $\N$. This shows that $e$ is
$G$-conjugate to $e_{\rm reg}$. Therefore, $M_j=\N(\z(e_{\rm
reg}))$ and $Z={\frak C}(e_{\rm reg})$, completing the proof. \qed

\smallskip

\noindent {\bf Remark.} It is quite possible that $\C(\g)$ is
equidimensional for any connected reductive group $G$ and for any
$p$. If this is the case, then $\dim \C(\g)=\dim {\frak C}(e_{\rm
reg})$, by Proposition~4.3. In [\cite{Spr}, Sect.~2], Springer
essentially computed $\dim \z(e_{\rm reg})$ for any semisimple
group $G$ in {\it bad} characteristic. In view of
Proposition~4.3(2), this computation yields a close formula for
$\dim {\frak C}(e_{\rm reg})$.

Proposition~4.3 also shows that if $\C(\g)$ is equidimensional,
then $$\dim \z(e)-\dim Z_G(e)-MT(\z(e))\,\le \,\dim \z(e_{\rm
reg})-\text{rk}\,G-\dim \z(\g)$$ for any $e\in\N$. Thus proving
the equidimensionality of $\C(\g)$ in bad characteristic would
have important implications for the detailed study of the adjoint
action of $G$.

\end{document}